\documentclass[10pt]{article}

\usepackage{amsmath,amssymb,amsthm,mathtools}
\usepackage{bm}
\usepackage{mathrsfs}

\usepackage[a4paper,margin=1in]{geometry}
\usepackage{microtype}
\usepackage{setspace}

\usepackage{graphicx}
\usepackage{booktabs}
\usepackage{float}

\usepackage[numbers]{natbib}

\usepackage[colorlinks=true,
            linkcolor=black,
            citecolor=blue,
            urlcolor=blue]{hyperref}

\usepackage{enumitem}
\usepackage{xcolor}
\usepackage{array}

\theoremstyle{definition}

\theoremstyle{remark}

\title{Invading and receding travelling waves of the Fisher-KPP equation with a mass-conserving, moving boundary}

\author{
Georgia R. Weatherley, Adrianne L. Jenner, and Michael C. Dallaston
}

\date{
School of Mathematical Sciences, Queensland University of Technology (QUT), Brisbane, Australia \\
\texttt{georgia.weatherley@hdr.qut.edu.au}\\
\today
}

\begin{document}

\maketitle

\begin{abstract}
The Fisher-KPP equation is the canonical reaction-diffusion equation used in the study of invasive phenomena in mathematical biology. While moving boundaries have already been included to introduce biologically realistic sharp interfaces, these flux-based conditions assume loss of mass as a byproduct of boundary movement. We consider a mass-conserving boundary condition at a free boundary of the Fisher-KPP equation to allow interface movement without population consumption. We use a combination of phase plane analysis, perturbation analysis, and numerical simulation to show that the model supports both invading and receding travelling wave solutions with distinct wave speeds and front densities. Prescribing the free boundary velocity as a linear function of the front density, we show the existence of multiple stable and unstable travelling wave solutions, as well as receding population blow-up. Our results indicate that in regions where multiple steady states exist, those with larger wave speeds are stable. This is corroborated by considering the parameter regime associated with the uniform steady state solution and performing linear stability analysis. The framework we develop is readily extendable to more general boundary velocity functions, offering a means of studying increasingly complex and biologically relevant boundary dynamics. 
\end{abstract}

\section{Introduction}
Migration, proliferation, and death processes drive the spatial expansion and contraction of biological populations. Traditionally, mathematical modelling has focused solely on invasion, with applications including the spread of invasive species \cite{shigesada1995modeling,skellam1951random}, the spread of cell populations in the context of cell experiments \cite{cai2007multi,warne2019using,crossley2023traveling} and malignant growth \cite{swanson2003virtual,swanson2008mathematical}. In contrast, comparatively less attention has been given to studying populations that recede. Across numerous disease contexts, recession indicates a shift towards meaningful disease outcomes such as stabilisation or remission. For example, we may expect the reduction of a tumour mass as cancer cells respond to innate immunological pressures \cite{wang2024mri,brindle2008new}. Reduced lesion volume in multiple sclerosis patients is evidence of the decreased inflammation and increased tissue repair that is typical of a remissive period of the disease. \cite{rovira2013magnetic,patrikios2006remyelination}. Recession can also be an indicator of disease stage. For instance, the regression of lipid-rich atherosclerosis plaques is typically more feasible in early stages, prior to the development of calcification \cite{williams2008rapid}.

\begin{figure}[!h]
\begin{center}
\def\svgwidth{\textwidth}
 \input{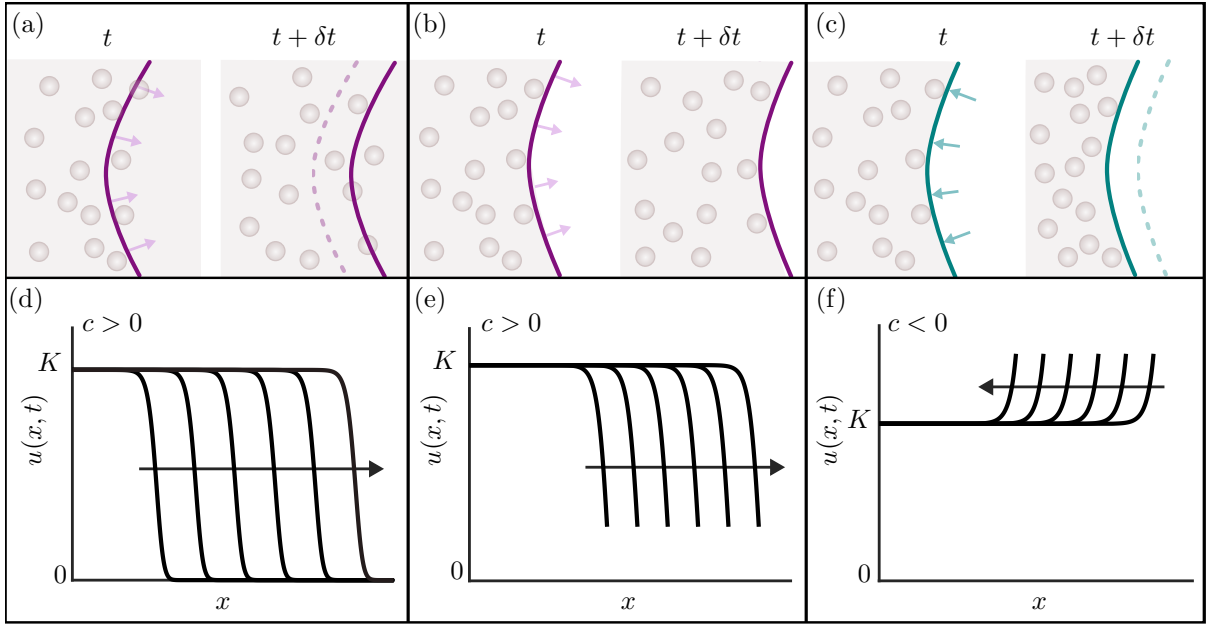}
    \caption{Schematic showing the evolution of a biological population that (a),(b) invades and (c) recedes, with the arrows indicating the direction of movement. (d-f) Travelling wave solutions corresponding to the density $u(x,t)$ of the populations in (a-c). In (a) the boundary of the population must be inferred by tracking an arbitrary threshold density, giving rise to smooth density fronts characteristic of the Fisher-KPP equation \eqref{F-KPP} in (d). (b),(c), represent the model \eqref{Our model dimen}-\eqref{1.4}, where the boundary of the domain changes based on the local population density \eqref{1.4} and mass is conserved \eqref{1.2}. This gives rise to the sharp, well-defined density fronts in (e),(f) respectively. The density of (f) increases above the carrying capacity $K$ due to the crowding effects introduced by the free boundary formulation \eqref{1.2}.}
    \label{fig: intro fig}
\end{center}
\end{figure}

 By representing a biological population as a moving density front, we are able to study its propensity to invade or recede (Fig. \ref{fig: intro fig}). Among the numerous systems of reaction-diffusion partial differential equations (PDEs) that have been developed to study the invasion of biological populations, the work of Fisher \cite{fisher1937wave} and Kolmogorov et al. \cite{kolmogorov1937study} remains canonical. Coupling linear diffusion with logistic growth, the Fisher–KPP equation can be written as
\begin{align}\label{F-KPP}
     \frac{\partial{u}}{\partial t} &= D\frac{\partial^2 u}{\partial x^2} + \lambda u\left(1-\frac{u}{K}\right),
\end{align}
where the population density $u(x,t)$ depends on position $x$ and time $t > 0$. The movement of individuals is described by linear diffusion with diffusivity $D > 0.$ The population grows logistically, with a proliferation rate $\lambda > 0$ and a carrying capacity density $K > 0.$ When considered on an infinite spatial domain ($-\infty < x < \infty$), a key property of the model is that a localised initial condition will evolve to a travelling wave solution with constant wave speed $c \geq 2\sqrt{D\lambda}$ in the long-time limit, $t \to \infty$. These well studied invading travelling wave solutions are smooth and lack a clearly defined front position (see \cite{Edelstein-Keshet_2005} and Fig. \ref{fig: intro fig}d). 

 While the Fisher-KPP model has proven itself to be a reliable tool for modelling invasive populations, key limitations of the model are that the population front is not well defined (Fig. \ref{fig: intro fig}a), and it cannot be used to study population recession. In fact, solutions of the Fisher-KPP model where $c < 2\sqrt{D\lambda}$ are discarded for introducing non-physical, negative population densities. This limitation has been overcome by extending the Fisher-KPP equation \eqref{F-KPP} to include a moving boundary, giving rise to the Fisher-Stefan model \cite{El-Hachem_McCue_Jin_Du_Simpson_2019,mccue2021exact,El-Hachem_McCue_Simpson_2021_I_R, McCue_2022, bui2024stability} (see Fig. \ref{fig: intro fig}b,c). Travelling wave solutions of the Fisher-Stefan model have both negative and positive wave speeds under well-defined, non-negative population densities. This formulation has also been modified to prescribe the front density as some constant $u_{f}$ \cite{El-Hachem_McCue_Simpson_2022_uf}. Examples of sharp-fronted solutions are shown in Fig. \ref{fig: intro fig}e,f. 

Governed by the Fisher-KPP equation, the Fisher-Stefan model describes the velocity of the free boundary at $x = L(t)$ by a Stefan condition (see Table \ref{Summary of BCs across models}). This condition was originally formulated to capture the evolving interface of two phases in the classical Stefan problem from thermodynamics \cite{hill1987one}. In the context of modelling a single biological species, it relates the velocity of the boundary to the flux of the population through the free boundary, leading to a loss of mass at the free boundary. In a biological context this could represent the population front interacting with a harsh environment, or the local conversion of the primary cell population to another unmodeled cell type at the front \cite{McCue_2022}. Given this condition is not appropriate in situations where the boundary is modified without loss of mass, 
one approach to retaining its physical significance is to incorporate a second biological population \cite{McCue_2022}. For instance, a primary population of cells may invade a background cell population with a moving interface separating the two populations, as examined by El-Hachem et al. \cite{el2020sharp}.

Mass-conserving boundary conditions arise naturally in models with both physical \cite{illingworth2005numerical,pegler2016dynamics} and biological applications \cite{gaffney1999modelling,chen2000free,baker2019free,murphy2021travelling}. In the work of Illingworth \& Golosnoy \cite{illingworth2005numerical}, interface movement counterbalances the flux of two diffusive solutes across the interface. To study the healing of corneal epithelial wounds, Gaffney et al. \cite{gaffney1999modelling} model a wound edge as a moving boundary, motivated by cell migration at the wound edge being influenced locally by the presence of a small, electric field. At the moving boundary, the cell density is fixed and the flux of cells into the boundary represents new tissue, advancing the position of the boundary while conserving local mass. In the mechanobiological model of epithelial tissue growth of Murphy et al. \cite{murphy2021travelling}, movement of a free boundary is coupled with underlying cell mechanics so that the tissue boundary expands without loss of material, and the cell density at the boundary is determined by the balance of mechanical forces. 

In this work, we develop a free-boundary model that captures the invasion and recession of a single biological population without assuming loss of mass. Akin to the Fisher-Stefan model, we reformulate the classical Fisher-KPP equation to include a moving boundary: 
\begin{align}\label{Our model dimen}
\frac{\partial{u}}{\partial t} &= D\frac{\partial^2 u}{\partial x^2} + \lambda u\left(1-\frac{u}{K}\right), \quad \text{on} \quad -\infty <x < L(t).
\end{align}
We take the moving boundary $x = L(t)$ to evolve according to a general function $f$ of cell density at this position:
\begin{align}\label{1.2}
\frac{\mathrm{d} L}{\mathrm{d}t} &= \left.{f}(u)\right|_{x=L(t)}.
\end{align}
This function may include the effects of cell degradation or deposition that change the boundary location, as well as external effects such as innate recovery of the domain. 
We assume a change in the position of the boundary is matched by the flux into the boundary, that is,
\begin{align}
    \frac{\partial u}{\partial x}(L(t),t) = -u(L(t),t)\frac{\mathrm{d}L}{\mathrm{d}t}.\label{1.4}
\end{align}
This condition ensures that loss of mass is not a byproduct of boundary movement, so that any changes in mass are solely attributed to the reaction term of the governing equation \eqref{Our model dimen}.

In this article, we examine the behaviour of a moving boundary Fisher-KPP equation \eqref{Our model dimen} under the mass-conserving condition \eqref{1.4} and a general boundary velocity function \eqref{1.2}, focusing on the presence and stability of travelling wave solutions. We use phase plane analysis to identify travelling wave solutions to this system in Section \ref{sec: phase plane analysis}, for which the boundary velocity function $f(u)$ does not need to be specified. We identify solutions with both positive, $c > 0$, and negative, $c < 0$, wave speeds to represent both invasion and recession (see Fig. \ref{fig: intro fig}). In Section \ref{sec: fu_formulation} we consider a linear boundary velocity function $f$ \eqref{1.2} under which invading and receding solutions can coexist. Our findings are supported by perturbation analysis, as well as numerical simulations of the PDE model in Section \ref{sec: PDE sims}, where we observe both stable and unstable travelling wave solutions. Finally, in Section \ref{sec: stability analysis} we analyse the stability of the steady state solution that exists when the population is at the carrying capacity; this analysis is representative of the stability properties of travelling waves that are observed for more general parameters in the numerical solutions.

\section{Phase plane analysis of travelling waves}\label{sec: phase plane analysis}
\subsection{Nondimensionalisation}
We simplify our analysis by rescaling the problem as follows: $x  \mapsto x\sqrt{\lambda/D}$, $t \mapsto \lambda t$, $u \mapsto u/K$, $L(t) \mapsto L(t)\sqrt{\lambda/D}$, and $f \mapsto \sqrt{D/\lambda }f$. The system (\ref{Our model dimen})-(\ref{1.4}) is then
\begin{align}
    &\frac{\partial{u}}{\partial t} = \frac{\partial^2 u}{\partial x^2} + u(1-u), \quad \text{on} \quad -\infty <x < L(t),
     \label{gov PDE}\\
      &\frac{\partial u}{\partial x}(L(t),t) = -u(L(t),t)\frac{\mathrm{d}L}{\mathrm{d}t},
     \label{BC right} \\
      &\frac{\mathrm{d}L}{\mathrm{d}t} = \left.{f}(u)\right|_{x=L(t)}.
    \label{BC L}  
\end{align} 
\subsection{Formulation of travelling wave solutions}
The governing equation \eqref{gov PDE} has been analysed extensively for travelling wave solutions to both the Fisher-KPP \cite{Canosa_1973,Murray_2002,Edelstein-Keshet_2005} and Fisher-Stefan \cite{El-Hachem_McCue_Jin_Du_Simpson_2019,El-Hachem_McCue_Simpson_2021_I_R, McCue_2022} models. 
To understand the effect of introducing a mass-conserving boundary condition \eqref{BC right} at $x = L(t)$, we start by undertaking a phase plane analysis to identify travelling wave solutions. We introduce a travelling wave coordinate $z = x-ct$, and let $L(t) = ct$, to seek solutions of the form $u(x,t) = U(z)$ to the governing equation \eqref{gov PDE}:
\begin{align}\label{2.6}
    \frac{\text{d}^2{U}}{\text{d} z^2} +c\frac{\text{d} U}{\text{d} z} + U(1-U) = 0, \quad \text{on} \quad -\infty < z < 0.
\end{align}
This second-order ODE \eqref{2.6} can be rewritten as a system of two autonomous, first-order ODEs,
\begin{align}\label{ODE for U}
    \frac{\text{d}{U}}{\text{d} z} &= W, \\
    \frac{\text{d}{W}}{\text{d} z} &= -cW - U(1-U).\label{ODE for W}  
\end{align}
Here, $z = 0$ is taken to be the location of the moving boundary without loss of generality. 
The appropriate boundary conditions are
\begin{align}
    &U(-\infty) = 1, \label{2.7} \\
     &W(0)= -U(0)c, \label{2.8}\\ 
    &c = \left.{f}(U)\right|_{z=0}. \label{2.9}
\end{align}
 The equilibrium points of the system \eqref{ODE for U}-\eqref{ODE for W} are $ (U,W) = (0,0)$ and $(U,W) = (1,0)$, with $(1,0)$ being a saddle point for all values of $c$. Conversely, the stability of $(0,0)$ is dependent on the wave speed $c$. It is a stable node if $c \geq 2$, a stable spiral if $0 < c < 2$, a centre if $c = 0$, an unstable spiral if $-2 < c < 0$, and an unstable node if $c \leq -2$ \cite{El-Hachem_McCue_Simpson_2021_I_R}.

\setlength{\extrarowheight}{10pt}
\begin{table}[!h]
\centering
\caption{Summary of boundary conditions at the free boundary $x = L(t)$ for the Fisher-Stefan model, the modified Fisher-Stefan model of El-Hachem et al. \cite{El-Hachem_McCue_Simpson_2022_uf}, and our model. We compare the relevant boundary conditions of the nondimensionalised PDE systems, their equivalent in the phase plane for equations \eqref{ODE for U}-\eqref{ODE for W}, as well as the truncation point of the travelling wave trajectory in the phase plane given by the respective boundary conditions.} 
\label{Summary of BCs across models}
\begin{tabular}{|p{4.6cm} p{2.8cm} p{2cm}|}
\hline
\textbf{PDE BCs} & \textbf{ODE BCs} & \textbf{Truncation} \\ \hline
\multicolumn{3}{|l|}{\textbf{Fisher-Stefan \cite{El-Hachem_McCue_Jin_Du_Simpson_2019,mccue2021exact,El-Hachem_McCue_Simpson_2021_I_R, McCue_2022, bui2024stability}:}} \\ 
$\displaystyle{u(L(t),t) = 0}$ & $U(0) = 0$ & $(0,-c/\kappa)$ \\
$\displaystyle{\frac{\mathrm{d}L}{\mathrm{d}t} = -\kappa\left.\frac{\partial u}{\partial x}\right|_{x=L(t)}}$ & $c = -\kappa W(0)$ &  \\ \hline
\multicolumn{3}{|l|}{\textbf{Modified Fisher-Stefan \cite{El-Hachem_McCue_Simpson_2022_uf}:}} \\ 
$\displaystyle{u(L(t),t) = u_{f}}$ & $U(0) = u_{f}$ & $(u_{f},-c/\kappa)$\\
$\displaystyle{\frac{\mathrm{d}L}{\mathrm{d}t} = -\kappa\left.\frac{\partial u}{\partial x}\right|_{x=L(t)}}$ &  $c = -\kappa W(0)$ &  \\ \hline
\multicolumn{3}{|l|}{\textbf{Mass-conserving:}} \\ 
 $\displaystyle{\frac{\partial u}{\partial x}(L(t),t) = -u(L(t),t)\frac{\mathrm{d}L}{\mathrm{d}t}}$ & $W(0) = -cU(0)$ & $(U^{*},-cU^{*})$ \\ 
 $\displaystyle{\frac{\mathrm{d}L}{\mathrm{d}t} = \left.{f}(u)\right|_{x=L(t)}}$ & $c = f(U)$ \rule[-10pt]{0pt}{20pt} &  \\ 
\hline
\end{tabular}
\vspace*{-1pt}
\end{table}%%%End of the table

For systems governed by the Fisher–KPP equation (\ref{F-KPP}), travelling wave trajectories are identical, departing $(1,0)$ along its unstable eigenvector to satisfy the condition (\ref{2.7}). These trajectories differ only in where they terminate, due to the distinct boundary conditions imposed by each model. For the Fisher-KPP model, these trajectories are heteroclinic orbits connecting $(1,0)$ and $(0,0)$ for $c \geq 2$, given non-physical negative population densities arise when $c < 2$. The inclusion of a moving boundary instead requires valid travelling wave trajectories to depart $(1,0)$ to intersect a boundary condition on $U(0)$ at some front density $U^*$. For the Fisher-Stefan model, the relevant boundary condition is $U(0) = 0$, so that $U^* = 0$, and the truncation occurs along the $W(z)$ axis at the point $(0,-c/\kappa)$, where $\kappa$ is the Stefan constant (see Table \ref{Summary of BCs across models}). Truncations away from the $W(z)$ axis (correspnding to a non-zero front density) have been examined in El-Hachem et al. \cite{El-Hachem_McCue_Simpson_2022_uf}, who consider the boundary condition $u(L(t),t) = u_{f}$. In this framework, the trajectories are required to intersect the vertical line $U(0) = u_{f}$ so that $U^* = u_{f}$, with $u_{f} = 0$ restoring the original Fisher-Stefan model.

\begin{figure}[!h]
 \def\svgwidth{\textwidth}
 \input{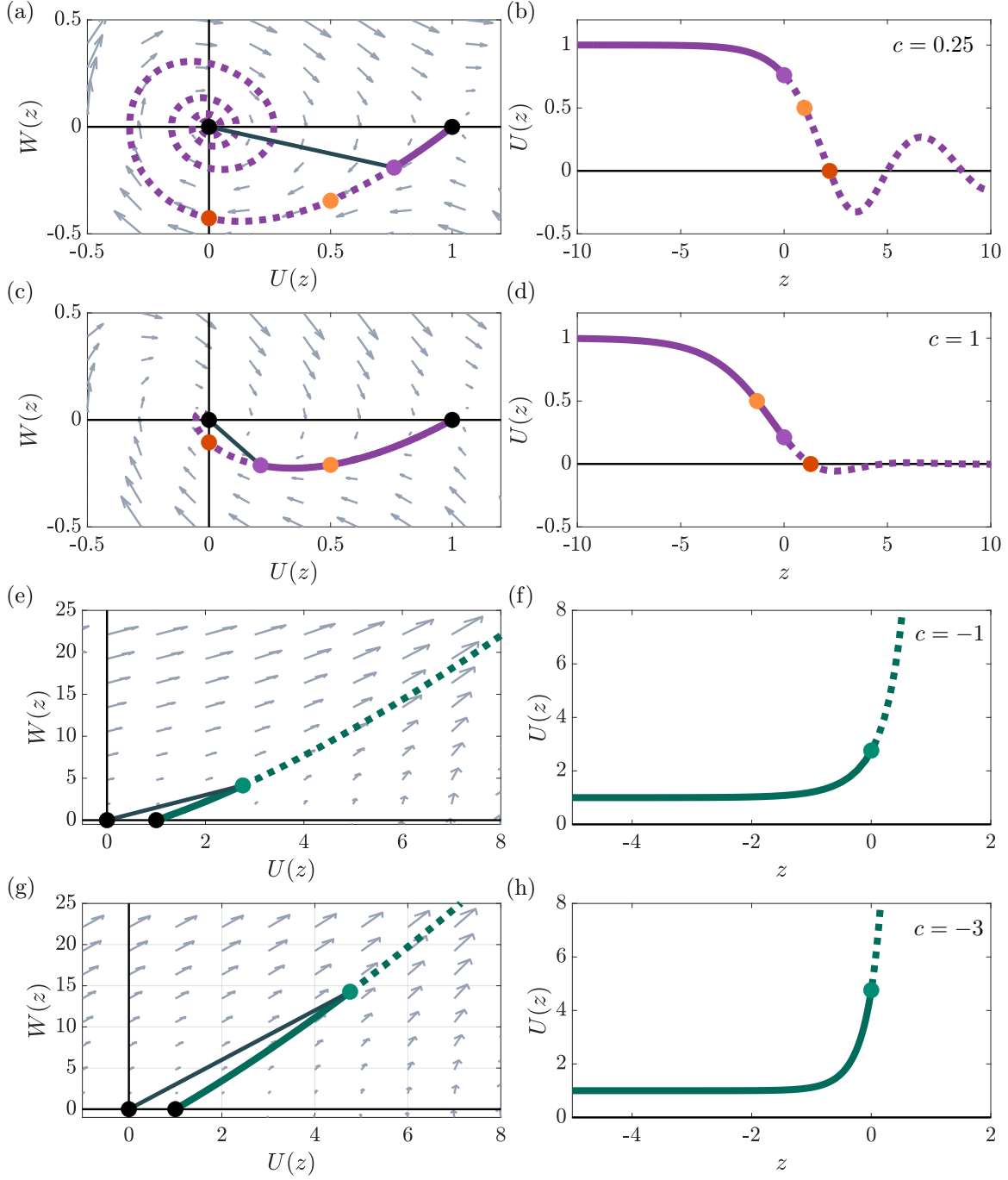}
    \caption{Travelling wave trajectories in $(U(z),W(z))$ phase space and their corresponding density profiles $U(z)$. Phase planes are shown in (a), (c), (e), and (g) for $c = 0.25, 1, -1,$ and $-3$, respectively, with corresponding density profiles in (b), (d), (f), and (h). The equilibrium points of the ODE system \eqref{ODE for U}–\eqref{ODE for W} are shown (black discs). Trajectories are coloured based on the sign of the wave speed (purple, green). Trajectories intersect the boundary condition $W = -cU$ \eqref{2.8} at points $(U^*,-cU^*)$ (purple, green coloured discs). The valid segment of each trajectory (solid line) corresponds to solutions of equations \eqref{gov PDE}–\eqref{BC L}, and beyond $(U^*,-cU^*)$ trajectories are truncated (dashed line). For positive wave speeds, additional truncation points are included: those associated with the boundary condition of El-Hachem et al. \cite{El-Hachem_McCue_Simpson_2022_uf} (yellow discs, where $U^* = 0.5$) and the Fisher–Stefan model \cite{El-Hachem_McCue_Jin_Du_Simpson_2019,El-Hachem_McCue_Simpson_2021_I_R, McCue_2022} (red discs, where $U^* = 0$). The vector field of the system is indicated (grey arrows).}
    \label{fig: result fig 1}
\end{figure}

In our model, the mass-conserving boundary condition \eqref{BC right} requires the trajectories to depart $(1,0)$ and intersect $W(z) = -c U(z)$ \eqref{2.8}. As a result, we truncate at the point $(U^{*},-c U^{*})$ in the $(U,W)$ plane. Unlike the Stefan condition which requires $U^* = 0$, the mass-conserving boundary condition introduces wave speed-dependent front density. To observe these truncated orbits in phase space, we initialise the ODE system \eqref{ODE for U}-\eqref{ODE for W} for a given $c$ near the unstable point $(1,0)$ and solve for $U(z)$ and $W(z)$. The mass-conserving condition (\ref{2.8}) allows both invading (see Fig. \ref{fig: result fig 1}a-d) and receding (see Fig. \ref{fig: result fig 1}e-h) travelling wave solutions. Importantly, these behaviours arise under strictly positive population densities. We confirm that for $c > 0$, travelling wave trajectories are truncations of the heteroclinic orbit between $(1,0)$ and $(0,0)$ (see Fig. \ref{fig: result fig 1}a,c). The case for $c < 0$ in Fig. \ref{fig: result fig 1}e,g shows that valid travelling wave trajectories depart $(1,0)$ and are truncated with positive front density $U^*$ to the right of $U(z) = 1$ in the $(U,W)$ phase space. 

\begin{figure}[!h]
\def\svgwidth{\textwidth}
 %% Creator: Inkscape 1.3 (0e150ed, 2023-07-21), www.inkscape.org
%% Accompanies image file 'XFigure_3.pdf' (pdf, eps, ps)
\begingroup%
  \makeatletter%
  \providecommand\color[2][]{%
    \errmessage{(Inkscape) Color is used for the text in Inkscape, but the package 'color.sty' is not loaded}%
    \renewcommand\color[2][]{}%
  }%
  \providecommand\transparent[1]{%
    \errmessage{(Inkscape) Transparency is used (non-zero) for the text in Inkscape, but the package 'transparent.sty' is not loaded}%
    \renewcommand\transparent[1]{}%
  }%
  \providecommand\rotatebox[2]{#2}%
  \newcommand*\fsize{\dimexpr\f@size pt\relax}%
  \newcommand*\lineheight[1]{\fontsize{\fsize}{#1\fsize}\selectfont}%
  \ifx\svgwidth\undefined%
    \setlength{\unitlength}{595.27559055bp}%
    \ifx\svgscale\undefined%
      \relax%
    \else%
      \setlength{\unitlength}{\unitlength * \real{\svgscale}}%
    \fi%
  \else%
    \setlength{\unitlength}{\svgwidth}%
  \fi%
  \global\let\svgwidth\undefined%
  \global\let\svgscale\undefined%
  \makeatother%
  \begin{picture}(1,0.88095238)%
    \lineheight{1}%
    \setlength\tabcolsep{0pt}%
    \put(0,0){\includegraphics[width=\unitlength,page=1]{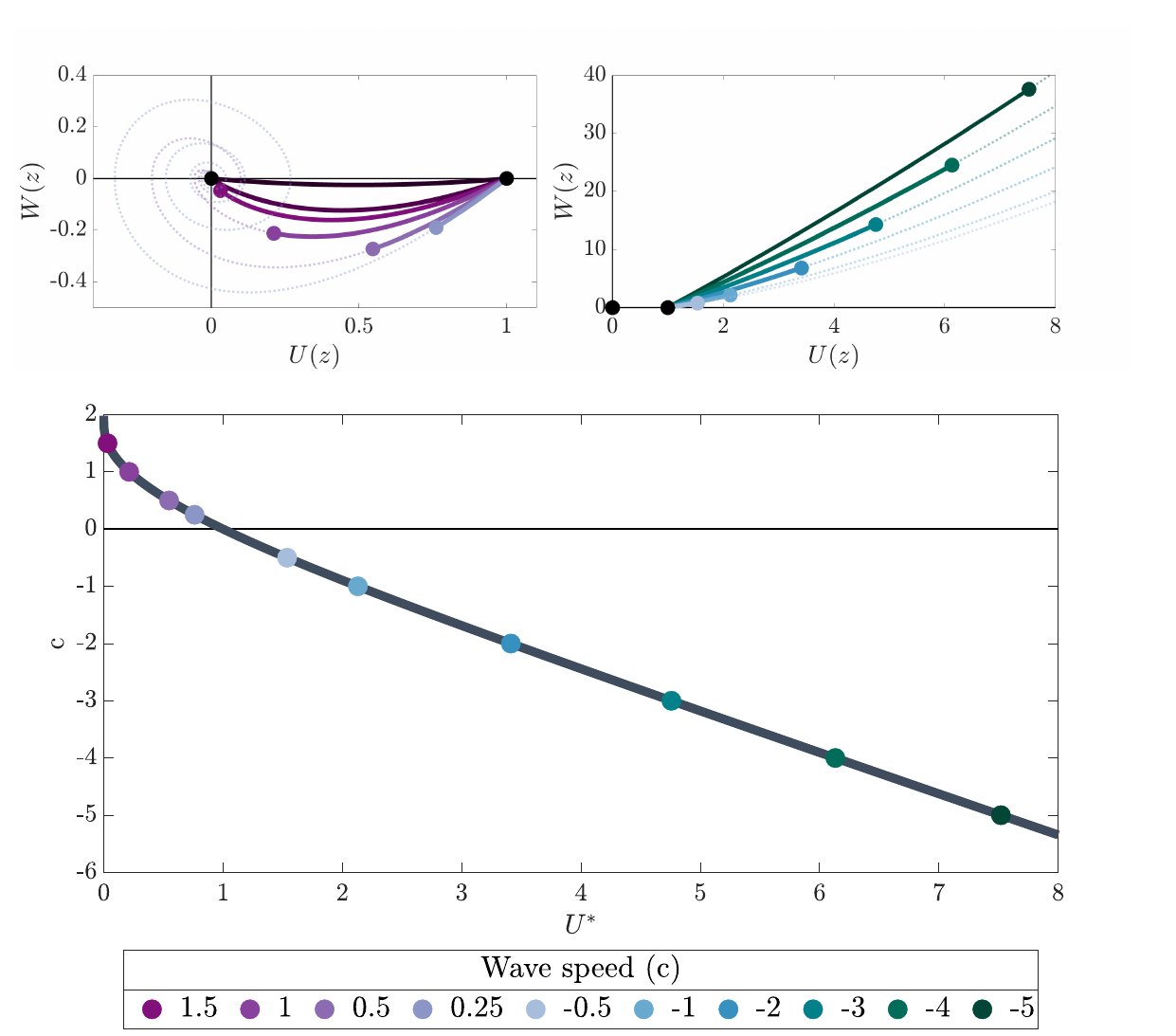}}%
    \put(0.02542601,0.8359085){\color[rgb]{0,0,0}\makebox(0,0)[t]{\lineheight{0.5}\smash{\begin{tabular}[t]{c}(a)\end{tabular}}}}%
    \put(0.03448165,0.54334957){\color[rgb]{0,0,0}\makebox(0,0)[t]{\lineheight{0.5}\smash{\begin{tabular}[t]{c}(c)\end{tabular}}}}%
    \put(0.47397232,0.8359085){\color[rgb]{0,0,0}\makebox(0,0)[t]{\lineheight{0.5}\smash{\begin{tabular}[t]{c}(b)\end{tabular}}}}%
  \end{picture}%
\endgroup%

    \caption{Family of trajectories for (a) positive and (b) negative wave speeds. Equilibrium points are shown (black discs), as well as intersections of trajectories with the boundary condition $W=-cU$ from equation \eqref{2.8} (coloured discs). Valid solutions to equations \eqref{gov PDE}–\eqref{BC L} run from $(1,0)$ to these truncation points $(U^*,-cU^*)$ (solid coloured lines), beyond which trajectories are truncated (dashed lines). In (a), trajectories connect $(1,0)$ to $(0,0)$ for $0 < c < 2$ (shown for $c = 0.25, 0.5, 1, 1.5, 2, 10$). No valid intersection occurs for $c = 2$ or $c = 10$ due to the linearisation about $(0,0)$. In (b), trajectories extend from $(1,0)$ into the first quadrant for $c < 0$ (shown for $c = -0.5, -1, -2, -3, -4, -5$). (c) shows the relationship between the front density $U^{*}$ and wave speed $c$, with phase plane intersections from (a) and (b) marked (coloured discs). This confirms that valid wave speeds satisfy $-\infty < c < 2$.}
    \label{fig: result fig 2}
\end{figure}

We observe restrictions on the wave speed by considering a broader range of wave speeds in Fig. \ref{fig: result fig 2}a,b. While the receding wave speed $c < 0$ is unrestricted, invading travelling wave solutions are restricted to $0 < c < 2$. For $c \geq 2$, trajectories approach the fixed point $(0,0)$ along an eigenvector with slope less negative than $-c$, so the trajectory cannot intersect the line $W=-cU$. This invading wave speed restriction of $0 < c < 2$ is also exhibited in the Fisher-Stefan model \cite{El-Hachem_McCue_Jin_Du_Simpson_2019}. It is interesting that for both models, the inclusion of the moving boundary allows us to consider wave speeds that are usually neglected for the Fisher-KPP model, where the minimum wave speed is $c > 2$.

In the density profiles $U(z)$ of Fig. \ref{fig: result fig 1}b,d, we observe for the invading case where $c > 0$ that $U^* \in (0,1)$. There is an inverse relationship between $c$ and $U(z)$, so that populations with low front densities $U^*$ invade with higher wave speeds $c$.  For the receding case, a higher front density $U^*$ is associated with a negative wave speed $c$ of a higher magnitude. As observed in the density profiles of $U(z)$ of Fig. \ref{fig: result fig 1}f,h, the recession of the boundary coupled with the mass-conserving condition induces a bunching behaviour near the population front so that $U^* > 1$. The relationship between wave speed $c$ and front density $U^{*}$ is shown in Fig. \ref{fig: result fig 2}c. The relationship is determined by numerically identifying $U^{*}$ for wave speeds $-6 < c < 2$. While this wave speed curve is nonlinear for $0 < c < 2$, it becomes asymptotically linear as $c \to -\infty$. 

\subsection{Perturbation solution for $c \to -\infty$}
\begin{figure}[!h]
\def\svgwidth{\textwidth}
 \input{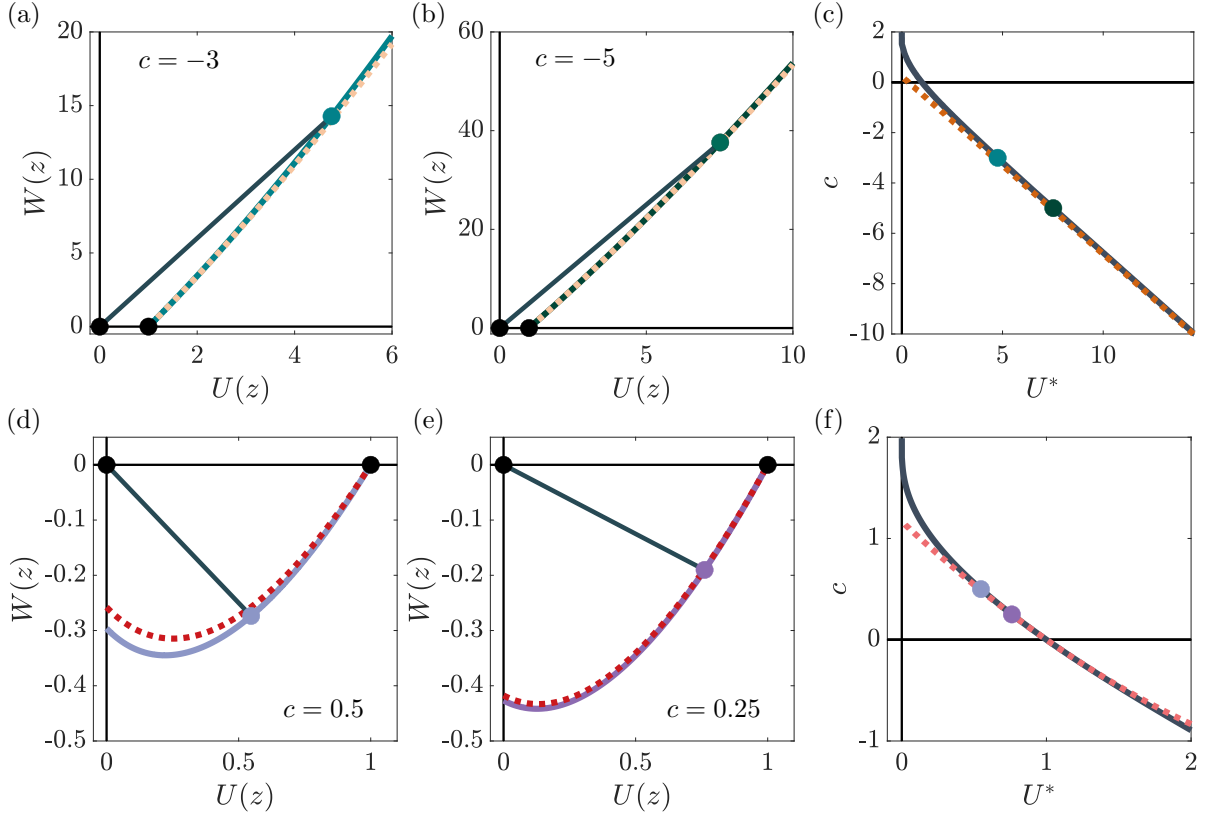}
    \caption{(a,b,c) Perturbation solutions for the phase plane and the relationship between $U^{*}$ and wave speed $c$ for $c \to -\infty$. In (a) and (b), phase plane trajectories are shown (solid green lines) for $c = -3$ and $c = -5$, alongside the $\mathcal{O}(c^{-3})$ perturbation solution (dashed yellow line) and the boundary condition from equation \eqref{2.8}. In (c), the $\mathcal{O}(1)$ perturbation solution is shown (dashed orange line) against the $U^{*}$ versus $c$ profile obtained via phase plane analysis (solid line). The truncation points from (a) and (b) are marked (coloured discs).(d,e,f) Perturbation solutions for the phase plane and the relationship between $U^{*}$ and wave speed $c$ for $|c| \ll 1$. In (d) and (e), phase plane trajectories are shown (solid purple lines) for $c = 0.5$ and $c = 0.25$, together with the $\mathcal{O}(c)$ perturbation solution (dashed red line) and the boundary condition from equation \eqref{2.8}. In (f), the $\mathcal{O}((U^*-1)^2)$ perturbation solution is shown (dashed orange line) against the $U^{*}$ versus $c$ profile obtained via phase plane analysis (solid line). The truncation points from (d) and (e) are indicated (coloured discs).}
    \label{fig: result fig 3 & 4}
\end{figure}

In the previous section, we numerically calculated the wave speed curve for the front density $U^{*}$ and wave speed $c$ (see Fig. \ref{fig: result fig 2}c). To gain analytic insights into this relationship we seek a perturbation solution for very negative wave speeds, $c \to -\infty$, by writing
\begin{align}
    W(U) = cW_{0} + W_{1} + \frac{1}{c}W_{2} + \frac{1}{c^2}W_{3} + \frac{1}{c^3}W_{4} + \mathcal{O}\left(\frac{1}{c^4}\right). \label{eq: pert form neg c}
\end{align}
The ODE system (\ref{ODE for U})-(\ref{ODE for W}) can be written as 
\begin{align}
  \frac{\text{d}{W}}{\text{d} U} = \frac{-cW - U(1-U)}{W}.\label{eq: DW DU first}
\end{align}
Substituting (\ref{eq: pert form neg c}) into \eqref{eq: DW DU first} and matching powers of $c$, we solve the system of differential equations for $W_{0}(U)$ through to $W_{4}(U)$ with the condition $W = 0$ when $U = 1$ to give
\begin{align*}
    W_{0}(U) &= 1 - U, \\
    W_{2}(U) &= \frac{-U^2}{2} + \frac{1}{2}, \\
    W_{4}(U) &= \frac{U^2}{4} + \frac{U^3}{6} - \frac{5}{12},
\end{align*}
with $W_{1}(U) = W_{3}(U) = 0$. We formulate $\mathcal{O}(c^{-1})$ and $\mathcal{O}(c^{-3})$ perturbation solutions by truncating equation (\ref{eq: pert form neg c}) at the desired order of $c$. To validate the accuracy of these solutions, we compare the $\mathcal{O}(c^{-3})$ perturbation solution against the numerically generated trajectory for $c = -3$ and $c = -5$ in the phase plane. As shown in Fig. \ref{fig: result fig 3 & 4}a,b, this perturbation solution accurately approximates the trajectory of the travelling wave solutions for these sufficiently negative wave speeds. 

To obtain an analytic expression for the wave speed curve, we equate the boundary condition $W = -cU$ \eqref{2.8} with the expansion for $W(U)$ that is accurate to $\mathcal{O}(c^{-3})$ to find 
\begin{align}
    -cU^{*} &= c(1-U^{*}) + \frac{1}{c}\left(\frac{-U^{*2}}{2} + \frac{1}{2}\right) + \frac{1}{c^3}\left(\frac{U^{*2}}{4} + \frac{U^{*3}}{6} - \frac{10}{24}\right) + \mathcal{O}\left(\frac{1}{c^{5}}\right).\label{eq: equating W}
\end{align}
Expanding $U^*$ in powers of $c$
\begin{align}
    U^{*} &= cU_{0}^{*} + U_{1}^{*}\label{eq: U expansion},
\end{align}
we substitute (\ref{eq: U expansion}) into (\ref{eq: equating W}) to find
\begin{align*}
    U_{0}^{*} &= -\sqrt2, \\
    U_{1}^{*} &= \frac{1}{3}.
\end{align*}
By (\ref{eq: U expansion}), we see that
\begin{align*}
    U^{*} &= -\sqrt2c + \frac{1}{3} + \mathcal{O}\left(\frac{1}{c}\right),
    \end{align*}
    which may be inverted to give
    \begin{align}
    c &= \frac{-U^{*}}{\sqrt2} + \frac{1}{3\sqrt2} + \mathcal{O}\left(\frac{1}{U^{*}}\right).\label{eq: pert slope sqrt2}
\end{align}
This approximation provides a linear estimate of the relationship between $c$ and $U^*$, which we show in Fig. \ref{fig: result fig 3 & 4}c against the wave speed curve determined numerically. The best agreement holds when the wave speed $c$ is very negative and the population has a high front density, $U^* \gg 1$. 

\subsection{Perturbation solution for $|c| \ll 1$}
Here, we build an analytic expression for the wave speed curve when $|c| \ll 1$. We have shown that the invading trajectories of our model match those of the Fisher-Stefan model, differing only in their point of truncation. We restate the perturbation solutions of (\ref{eq: DW DU first}) for $c \ll 1$,
\begin{align}
    W(U) = W_{0}(U) + cW_{1}(U) + \mathcal{O}(c^2) \label{eq: pert small c}.
\end{align}
Integrating the differential equations for $W_{0}$ and $W_{1}$  obtained by substituting equation (\ref{eq: pert small c}) into equation (\ref{eq: DW DU first}), and using the initial condition $U(1) = 0$, gives \cite{El-Hachem_McCue_Jin_Du_Simpson_2019, El-Hachem_McCue_Simpson_2021_I_R}:
\begin{align}
W_{0}(U) &=  -(1-U)\sqrt{\frac{1+2U}{3}}, \label{eq: for SS sol} 
\end{align}
\begin{align*}
W_{1}(U) &= \frac{-(U-2)(1+2U)^{3/2} - \sqrt{27}}{5(U-1)\sqrt{1 + 2U}}.
\end{align*}
We have both $\mathcal{O}(1)$ and $\mathcal{O}(c)$ perturbation solutions for small $c$. In Fig. \ref{fig: result fig 3 & 4}d,e, we see that the $\mathcal{O}(c)$ solution is accurate for $c = 0.25$ and $c = 0.5$ by comparing to numerically generated phase plane trajectories. The expansion for $W(U)$ accurate to $\mathcal{O}(1)$ \eqref{eq: for SS sol} is exact for $c = 0$. In Fig. \ref{fig: result fig 2}c, the steady state solution $c = 0$ has front density is $U^* = 1$. Recalling the mass-conserving condition $W = -cU$ \eqref{2.8}, we satisfy $c = 0$ by setting $W(U) = 0$ in \eqref{eq: for SS sol} to identify the root of \eqref{eq: for SS sol}, $U(z) = 1$, as the uniform, steady state solution.

An analytic expression for the wave speed curve for $|c| \ll 1$ is developed by expanding in powers of $c$
\begin{align*}
    U^{*} &= 1 + cU_{1}^{*} + c^2U_{2}^{*}.
\end{align*}
We equate $W = -cU$ with the expansion for $W(U)$ given by \eqref{eq: pert small c}, to identify
\begin{align*}
    U_{1}^{*} &= -1, \\
    U_{2}^{*} &= \frac{1}{6}.
\end{align*}
The expansion for $U^*$ is 
\begin{align*}
    U^{*} &= 1 - c + \frac{1}{6}c^2,
\end{align*}
which may be inverted to give
\begin{align} 
    c &= -(U^{*}-1) + \frac{1}{6}(U^{*}-1)^2 + \mathcal{O}((U^{*}-1)^3). \label{eq: small_c_pert}
\end{align}
In Fig. \ref{fig: result fig 3 & 4}f we show the expansion \eqref{eq: small_c_pert} against the wave speed curve determined numerically. It captures the nonlinear relationship observed between $U^*$ and $c$ in Fig. \ref{fig: result fig 2}c. As expected, the best agreement holds when the population invades or recedes slowly with $|c| \ll 1$ and a front density close to the carrying capacity $U = 1$.

\section{Travelling Waves for linear boundary velocity function}\label{sec: fu_formulation}
We have confirmed the existence of travelling wave solutions and have a clear understanding of the relationship between their wave speeds and front densities. This analysis is independent of our choice of boundary velocity function \eqref{BC L}. Recall $f(u)$ to be some general function of the population density at the interface that describes the motion of the boundary. In this work, we consider a simple linear formulation
\begin{align}
    \frac{\mathrm{d}L}{\mathrm{d}t} = \left.{f}(u)\right|_{x=L(t)} = \alpha u(L(t),t) + \beta,
    \label{BC_fu} 
\end{align}
where $\alpha$ and $\beta$ are nondimensional constants which may be positive or negative. The slope $\alpha$ may be interpreted as scaling the population density at the interface to capture its effects through degradation or deposition of material. The intrinsic velocity $\beta$ may represent external effects, such as a natural recovery of the host environment. For small front densities, the behaviour of the system is dominated by $\beta$ and the direction of movement will largely depend on whether external effects encourage recovery ($\beta < 0)$ or degradation $(\beta > 0)$. For large front densities, $\alpha$ will dominate and the direction of movement will largely depend on whether individuals at the boundary have a degenerative effect $(\alpha > 0)$ or reparative effect $(\alpha < 0)$. 

\begin{figure}[!h]
\def\svgwidth{\textwidth}
 %% Creator: Inkscape 1.3 (0e150ed, 2023-07-21), www.inkscape.org
%% Accompanies image file 'XFigure_5.pdf' (pdf, eps, ps)
\begingroup%
  \makeatletter%
  \providecommand\color[2][]{%
    \errmessage{(Inkscape) Color is used for the text in Inkscape, but the package 'color.sty' is not loaded}%
    \renewcommand\color[2][]{}%
  }%
  \providecommand\transparent[1]{%
    \errmessage{(Inkscape) Transparency is used (non-zero) for the text in Inkscape, but the package 'transparent.sty' is not loaded}%
    \renewcommand\transparent[1]{}%
  }%
  \providecommand\rotatebox[2]{#2}%
  \newcommand*\fsize{\dimexpr\f@size pt\relax}%
  \newcommand*\lineheight[1]{\fontsize{\fsize}{#1\fsize}\selectfont}%
  \ifx\svgwidth\undefined%
    \setlength{\unitlength}{593.85506698bp}%
    \ifx\svgscale\undefined%
      \relax%
    \else%
      \setlength{\unitlength}{\unitlength * \real{\svgscale}}%
    \fi%
  \else%
    \setlength{\unitlength}{\svgwidth}%
  \fi%
  \global\let\svgwidth\undefined%
  \global\let\svgscale\undefined%
  \makeatother%
  \begin{picture}(1,0.37951541)%
    \lineheight{1}%
    \setlength\tabcolsep{0pt}%
    \put(0,0){\includegraphics[width=\unitlength,page=1]{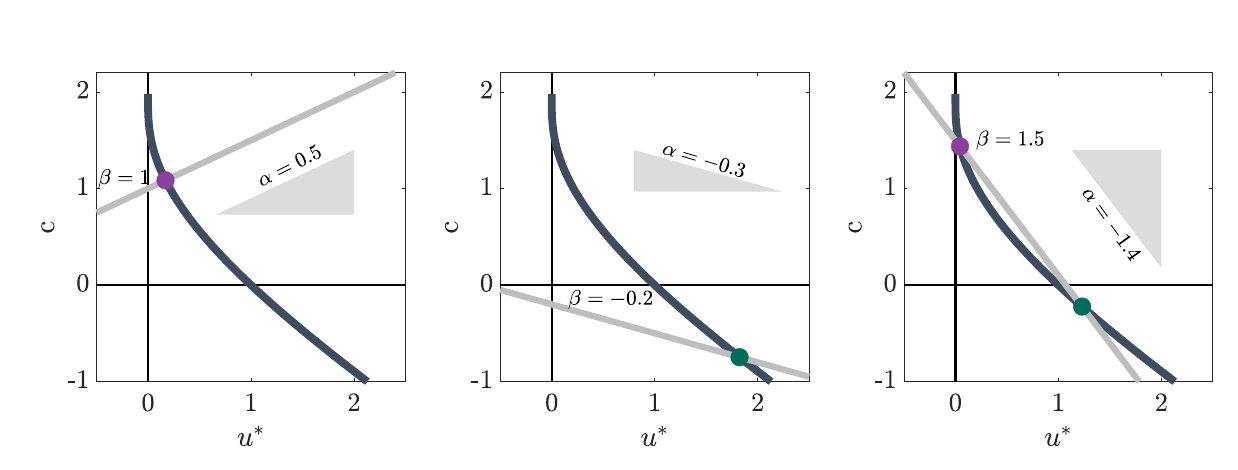}}%
    \put(0.68726217,0.3474737){\color[rgb]{0,0,0}\makebox(0,0)[t]{\lineheight{0.5}\smash{\begin{tabular}[t]{c}(c)\end{tabular}}}}%
    \put(0.03176395,0.34498836){\color[rgb]{0,0,0}\makebox(0,0)[t]{\lineheight{0.5}\smash{\begin{tabular}[t]{c}(a)\end{tabular}}}}%
    \put(0.35745187,0.34773479){\color[rgb]{0,0,0}\makebox(0,0)[t]{\lineheight{0.5}\smash{\begin{tabular}[t]{c}(b)\end{tabular}}}}%
  \end{picture}%
\endgroup%

    \caption{Intersection point $(u^*,c)$ of boundary velocity function (\ref{BC_fu}), with the wave speed curve (Fig. \ref{fig: result fig 2}c) for different choices of $\alpha$ and $\beta$. (a) $\alpha = 0.5$ and $\beta = 1$ with intersection point $(0.17,1.08)$. (b) $\alpha = -0.3$ and $\beta =-0.2$ with intersection point $(1.82,-0.75)$. (c) $\alpha = -1.4$ and $\beta = 1.5$ with two intersection points: $(1.23,-0.22)$ and $(0.04,1.44)$. Each plot includes a line segment of slope $\alpha$. Intersection points correspond to travelling wave solutions and are coloured by the sign of $c$, where purple indicates $c > 0$ and green indicates $c < 0$.}
    \label{fig: BC intersections}
\end{figure}
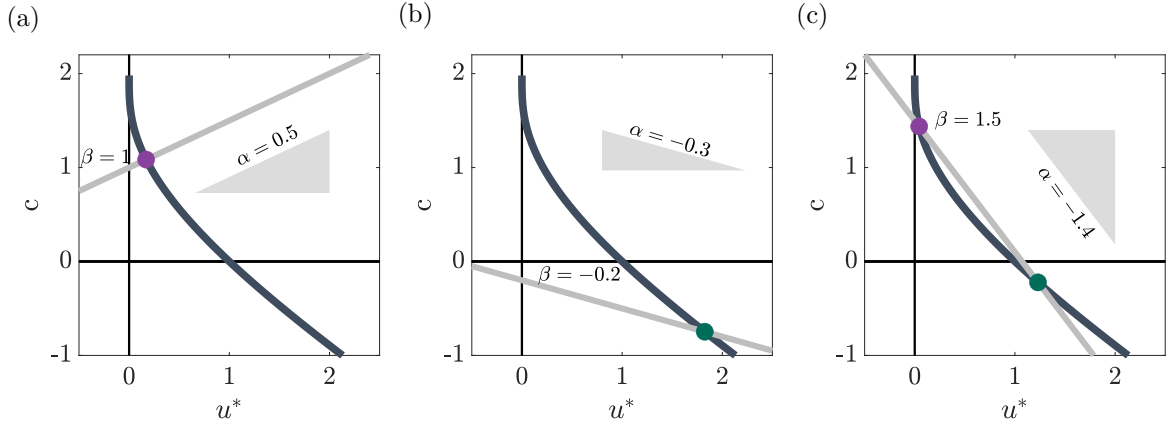

For a given boundary velocity function $f$ (\ref{BC_fu}), relevant travelling wave solutions exist at points where the boundary velocity function \eqref{BC_fu} intersects the wave speed curve depicted in Fig. \ref{fig: result fig 2}c. For a given $(\alpha$, $\beta)$ pair, the intersection point $(u^{*}, c)$ characterises the front density and wave speed of the travelling wave solution to the PDE system \eqref{gov PDE}-\eqref{BC L}. Multiple intersections, each corresponding to a travelling wave solution, are possible due to the concave-up relationship between $u^{*}$ and $c$. The values of $\alpha$ and $\beta$, determine whether the linear boundary velocity function intersects the wave speed curve zero, one, or two times. We represent intersections of single invading and receding solutions in Fig. \ref{fig: BC intersections}a and Fig. \ref{fig: BC intersections}b respectively, showing in Fig. \ref{fig: BC intersections}c that, among other possibilities, a single $(\alpha$, $\beta)$ pair may be associated with both invading and receding behaviours.

\begin{figure}[!h]
\def\svgwidth{\textwidth}
 \input{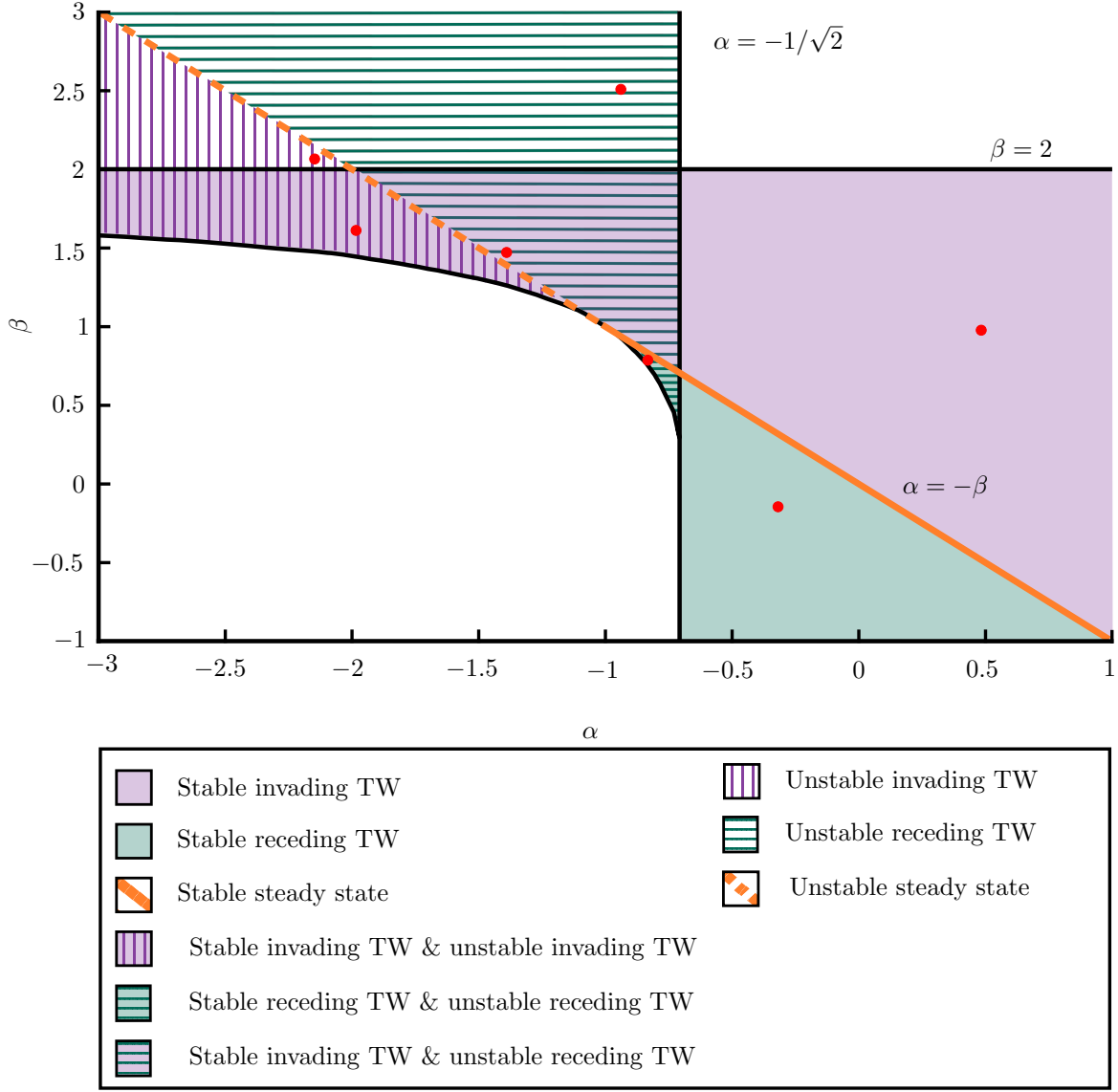}
    \caption{Travelling wave existence and stability for parameters $\alpha$ and $\beta$ of the boundary velocity function (\ref{BC_fu}). We consider ranges $-3 < \alpha < 1$ and $-1 < \beta < 3$. Intersections between the boundary velocity function and wave speed curve determine the existence of travelling wave solutions. We recover regions with no valid solution (white), invading travelling wave solutions (purple), receding travelling wave solutions (green), and steady state solutions (orange). We highlight (red discs) the three cases demonstrated in Fig. \ref{fig: BC intersections}, as well as those simulated in subsequent numerical simulations of the PDE system (see Figs. \ref{fig: PDE TW}-\ref{fig: unstable PDE}). Solid colouring represents a stable solution, while dashed colouring represents an unstable solution. In some regions we see that multiple solutions are possible, but only one is stable. Boundaries (black lines) between various regions that are known analytically from the perturbation solutions shown in Fig. \ref{fig: result fig 3 & 4} are labelled, with the remaining boundaries determined numerically.}
    \label{fig: paramspace fig}
\end{figure}
 By repeating the process demonstrated in Fig. \ref{fig: BC intersections} across different 
$(\alpha$, $\beta)$ pairs, we determine (numerically) the relationship between $\alpha$ and $\beta$ and the number and speed of travelling wave solutions. In Fig.
\ref{fig: paramspace fig} we summarise these results in the $(\alpha, \beta)$ plane. A large number of $(\alpha,\beta)$ pairs will give rise to no or a single travelling wave solution. For particular finite regions in the $(\alpha, \beta)$ plane, there are multiple travelling wave solutions, which can include two advancing travelling waves, two receding travelling waves, and one of each (Fig. \ref{fig: BC intersections}c). 

The boundaries between the different regions of Fig.
\ref{fig: paramspace fig} can mostly be determined analytically. The line $\alpha = -\beta$ corresponds to the steady state $(u^*,c) = (1,0)$ and marks the transition between invading ($c>0$) and receding ($c<0$) waves. The line $\beta = 2$ arises from the upper bound on the wave speed $c < 2$, above which intersections require sufficiently negative slopes $\alpha$. When $\beta > 2$, travelling wave solutions exist only if $\alpha < -1/\sqrt{2}$, reflecting the asymptotic slope of the wave speed curve for $c \to - \infty$ \eqref{eq: pert slope sqrt2}. As the relationship between $u^{*}$ and $c$ is independent of the boundary velocity function $f(u)$, we can numerically simulate the PDE system \eqref{gov PDE}-\eqref{BC L} under informed parameterisations of the boundary velocity function $f(u)$ (\ref{BC_fu}) based on the behaviours we want the system to exhibit. 

\section{PDE solutions}\label{sec: PDE sims}
In this section, we present numerical solutions of the PDE system \eqref{gov PDE}–\eqref{BC L}. As established in the previous section, a linear boundary velocity function (\ref{BC_fu}) gives rise to single travelling wave, multiple travelling wave, and non-travelling wave behaviours. To observe these different regimes, we impose conditions on $\alpha$ and $\beta$, noting from Fig. \ref{fig: paramspace fig} the significance of $\alpha = -1/\sqrt{2}$ and $\beta = 2$ in defining the transition between single and multiple travelling wave solutions. 

\subsection{Front-fixing and numerics}
To simulate the system \eqref{gov PDE}-\eqref{BC L} numerically, we solve over the domain $0 < x < L(t)$ and impose a Neumann boundary condition at $x = 0$. We begin by introducing a Landau-type transformation $\rho = x/L(t)$ \cite{landau1950heat,Back_McCue_Moroney_2014,illingworth2005numerical}. This transformation maps the moving boundary problem of domain $0 < x < L(t)$ to the fixed domain $0 < \rho < 1$. Defining $u(x,t) = \mathcal{U}(\rho(x,t),t)$, it follows that 
\begin{align}
    \frac{\partial{\mathcal{U}}}{\partial t} = \frac{1}{L^2}\frac{\partial^2 \mathcal{U}}{\partial \rho^2} +  \frac{\rho}{L}\frac{\text{d}L}{\text{d}t} \frac{\partial \mathcal{U}}{\partial \rho}+ \mathcal{U}(1-\mathcal{U}), \quad \text{on} \quad 0 <\rho < 1, \label{S1.1}    
\end{align}
subject to boundary conditions,
\begin{align}
    \frac{\partial \mathcal{U}}{\partial \rho}(0,t) &= 0 \label{S1.2},\\ 
    \frac{\partial \mathcal{U}}{\partial \rho}(1,t) &= -\mathcal{U}L\frac{\mathrm{d}L}{\mathrm{d}t}, \label{S1.3}\\ 
    \frac{\mathrm{d}L}{\mathrm{d}t} &= \left.f(\mathcal{U})\right|_{\rho=1} \label{S1.4}.
\end{align} 
To numerically solve \eqref{S1.1}-\eqref{S1.4}, we use a finite difference method for spatial discretisation, and then MATLAB's ode15s to integrate the resulting ODE system through time. We discretise the interval $0 < \rho < 1$ into $N$ node points of uniform spacing $\delta \rho = 1/(N-1)$. The $i$th mesh point is given by $\rho_{i} = (i-1)\delta \rho$ for $i = 1, ..., N$, so that at time $t$, the density $\mathcal{U}(\rho _{i},t)$ is represented by $\mathcal{U}_{i}$. To handle the spatial derivatives we use central difference approximations
\begin{align*}
    \frac{\partial \mathcal{U}_{i}}{\partial \rho} &\approx \frac{\mathcal{U}_{i+1} - \mathcal{U}_{i-1}}{2\delta\rho}, \\
    \frac{\partial^2 \mathcal{U}_{i}}{\partial \rho^2}&\approx \frac{\mathcal{U}_{i+1} - 2\mathcal{U}_{i} + \mathcal{U}_{i-1}}{\delta \rho^2},
\end{align*}
to give 
\begin{align*}
\frac{\text{d}{\mathcal{U}}_{i}}{\text{d}t} = 
 \begin{cases} 
      \displaystyle{\frac{1}{L^2}\frac{2\mathcal{U}_{2} -2\mathcal{U}_{1}}{\delta \rho^2}  + \mathcal{U}_{1}(1- \mathcal{U}_{1})} \quad \quad i = 1,& \\
      \displaystyle{\frac{1}{L^2}\frac{\mathcal{U}_{i+1} - 2\mathcal{U}_{i} + \mathcal{U}_{i-1}}{\delta \rho^2} + \frac{\mathrm{d}L}{\mathrm{d}t}\frac{\rho}{L}\frac{\mathcal{U}_{i+1} - \mathcal{U}_{i-1}}{2\delta\rho} + \mathcal{U}_{i}(1- \mathcal{U}_{i})} \quad \quad 2 \leq i \leq N-1,& \\
      \displaystyle{\frac{1}{L^2}\frac{\mathcal{U}_{N-1} - 2\delta \rho \mathcal{U}_{N}L\frac{\mathrm{d}L}{\mathrm{d}t} -2\mathcal{U}_{N} + \mathcal{U}_{N-1}}{\delta \rho^2} + \frac{\mathrm{d}L}{\mathrm{d}t}\frac{\rho}{L}\left(-\mathcal{U}_{N}L \frac{\mathrm{d}L}{\mathrm{d}t}\right) + \mathcal{U}_{N}(1- \mathcal{U}_{N})} \quad \quad i = N.& \\
   \end{cases}
\end{align*} 
Here, $\mathrm{d}L/\mathrm{d}t$ refers the boundary velocity function, which is given by 
\begin{equation*}
    \frac{\mathrm{d}L}{\mathrm{d}t} = \frac{\mathrm{d}\mathcal{U}_{N+1}}{\mathrm{d}t} = f(U_{N}) = \alpha \mathcal{U}_{N} + \beta.
\end{equation*}
The MATLAB code is available on \href{https://github.com/Georgia-Weatherley/travelling_waves_mass_conserving_boundary_condition}{github}.

\subsection{Single travelling wave solutions}
To simulate a single travelling wave solution, we sample $(\alpha,\beta)$ pairs where $\alpha > -1/\sqrt2$ and $\beta < 2$ (see Fig. \ref{fig: paramspace fig}). Receding solutions require $\alpha < -\beta$, while invading solutions require  
$\alpha > -\beta$. We use a uniform initial condition as follows:
\begin{align*}
   u(x, 0) = 0.5, \quad 0<x<L(t). 
\end{align*}
To prevent the left-hand boundary at $x = 0$ from affecting travelling wave solutions, $L(0)$ is chosen depending on whether the invading or receding case is considered. For the invading case,  $L(0) = 25$ is sufficiently large for the wave to stabilise away from the fixed boundary at $x = 0$. For the receding case, $L(0) = 100$ allows the wave form to be sufficiently observed before it approaches the boundary at $x = 0$.
\begin{figure}[!h]
\def\svgwidth{\textwidth}
 %% Creator: Inkscape 1.3 (0e150ed, 2023-07-21), www.inkscape.org
%% Accompanies image file 'XFigure_7.pdf' (pdf, eps, ps)
\begingroup%
  \makeatletter%
  \providecommand\color[2][]{%
    \errmessage{(Inkscape) Color is used for the text in Inkscape, but the package 'color.sty' is not loaded}%
    \renewcommand\color[2][]{}%
  }%
  \providecommand\transparent[1]{%
    \errmessage{(Inkscape) Transparency is used (non-zero) for the text in Inkscape, but the package 'transparent.sty' is not loaded}%
    \renewcommand\transparent[1]{}%
  }%
  \providecommand\rotatebox[2]{#2}%
  \newcommand*\fsize{\dimexpr\f@size pt\relax}%
  \newcommand*\lineheight[1]{\fontsize{\fsize}{#1\fsize}\selectfont}%
  \ifx\svgwidth\undefined%
    \setlength{\unitlength}{562.64351432bp}%
    \ifx\svgscale\undefined%
      \relax%
    \else%
      \setlength{\unitlength}{\unitlength * \real{\svgscale}}%
    \fi%
  \else%
    \setlength{\unitlength}{\svgwidth}%
  \fi%
  \global\let\svgwidth\undefined%
  \global\let\svgscale\undefined%
  \makeatother%
  \begin{picture}(1,0.59689991)%
    \lineheight{1}%
    \setlength\tabcolsep{0pt}%
    \put(0,0){\includegraphics[width=\unitlength,page=1]{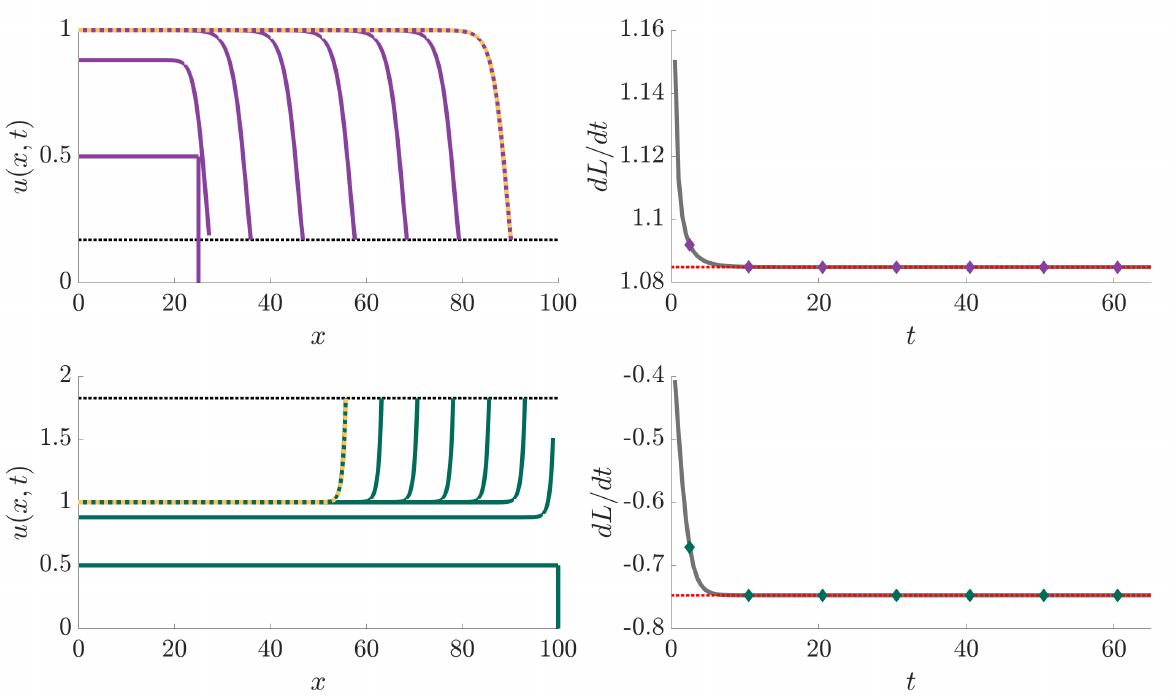}}%
    \put(0.01357074,0.56530207){\color[rgb]{0,0,0}\makebox(0,0)[t]{\lineheight{0.5}\smash{\begin{tabular}[t]{c}(a)\end{tabular}}}}%
    \put(0.01357074,0.29733158){\color[rgb]{0,0,0}\makebox(0,0)[t]{\lineheight{0.5}\smash{\begin{tabular}[t]{c}(c)\end{tabular}}}}%
    \put(0.49346372,0.56530207){\color[rgb]{0,0,0}\makebox(0,0)[t]{\lineheight{0.5}\smash{\begin{tabular}[t]{c}(b)\end{tabular}}}}%
    \put(0.49346372,0.29733158){\color[rgb]{0,0,0}\makebox(0,0)[t]{\lineheight{0.5}\smash{\begin{tabular}[t]{c}(d)\end{tabular}}}}%
    \put(1.312355,0.5434078){\color[rgb]{0,0,0}\makebox(0,0)[t]{\lineheight{0.5}\smash{\begin{tabular}[t]{c}$c = 0.25$\end{tabular}}}}%
    \put(1.32559028,0.44186959){\color[rgb]{0,0,0}\makebox(0,0)[t]{\lineheight{0.5}\smash{\begin{tabular}[t]{c}$c = 1$\end{tabular}}}}%
    \put(0,0){\includegraphics[width=\unitlength,page=2]{XFigure_7.pdf}}%
  \end{picture}%
\endgroup%

    \caption{(a,c) Numerical solutions of the PDE system \eqref{gov PDE}–\eqref{BC L} for different values of $\alpha$ and $\beta$, shown at $t = [0, 2, 10, 20, 30, 40, 50, 60]$. Direction of travel is marked (arrow). Superimposed are the wave profiles from the numerical phase plane trajectory (dashed yellow line) and the expected front density $u^*$ from phase plane analysis (dashed black line). (b,d) Front speeds at $x = L(t)$, with time points from (a,c) indicated (coloured markers). The expected wave speed $c$ from the phase plane analysis is shown (dashed red line). In (a,b), $\alpha = 0.5$, $\beta = 1$, giving rise to an invading travelling wave with $u^* = 0.17$ and $c = 1.08$. In (c,d), $\alpha = -0.3$, $\beta = -0.2$, giving rise to a receding travelling wave with $u^* = 1.82$ and $c = -0.75$. Colours (purple for $c > 0$, green for $c < 0$) indicate the sign of $c$. Note that these simulations correspond to the parameters shown in Fig. \ref{fig: BC intersections}a,b.}
    \label{fig: PDE TW}
\end{figure}

In Fig. \ref{fig: PDE TW} we show the evolution of the population density $u(x, t)$ for two different $(\alpha,\beta)$ pairs. In both instances, we expect the uniform initial condition to evolve into a travelling wave of constant speed. In Fig. \ref{fig: PDE TW}a,b we set $\alpha = 0.5$ and $\beta = 1$, to show a solution that approaches an invading travelling wave of speed $c = 1.08$ with front density $u^* = 0.17$. Similarly, in Fig. \ref{fig: PDE TW}c,d we set $\alpha = -0.3$ and $\beta = -0.2$ to give rise to a receding travelling wave of speed $c = -0.75$ with front density $u^* = 1.82$. By sampling the appropriate parameter spaces of Fig. \ref{fig: paramspace fig}, we have simulated two travelling wave solutions whose characteristics match those predicted by the theoretical phase plane solutions. As was observed in the density profiles of Fig. \ref{fig: result fig 1}, the front density of the invading solution in Fig. \ref{fig: PDE TW}a reduces below the steady state density of $u = 1$ so that $u^* \in (0,1).$ In contrast, for the receding case shown in Fig. \ref{fig: PDE TW}c, the front density increases above $u = 1$ so that  $u^* > 1.$

\subsection{Multiple travelling wave solutions}
\begin{figure}[!h]
\def\svgwidth{\textwidth}
 %% Creator: Inkscape 1.3 (0e150ed, 2023-07-21), www.inkscape.org
%% Accompanies image file 'XFigure_8.pdf' (pdf, eps, ps)
\begingroup%
  \makeatletter%
  \providecommand\color[2][]{%
    \errmessage{(Inkscape) Color is used for the text in Inkscape, but the package 'color.sty' is not loaded}%
    \renewcommand\color[2][]{}%
  }%
  \providecommand\transparent[1]{%
    \errmessage{(Inkscape) Transparency is used (non-zero) for the text in Inkscape, but the package 'transparent.sty' is not loaded}%
    \renewcommand\transparent[1]{}%
  }%
  \providecommand\rotatebox[2]{#2}%
  \newcommand*\fsize{\dimexpr\f@size pt\relax}%
  \newcommand*\lineheight[1]{\fontsize{\fsize}{#1\fsize}\selectfont}%
  \ifx\svgwidth\undefined%
    \setlength{\unitlength}{570.92871959bp}%
    \ifx\svgscale\undefined%
      \relax%
    \else%
      \setlength{\unitlength}{\unitlength * \real{\svgscale}}%
    \fi%
  \else%
    \setlength{\unitlength}{\svgwidth}%
  \fi%
  \global\let\svgwidth\undefined%
  \global\let\svgscale\undefined%
  \makeatother%
  \begin{picture}(1,0.85899184)%
    \lineheight{1}%
    \setlength\tabcolsep{0pt}%
    \put(0,0){\includegraphics[width=\unitlength,page=1]{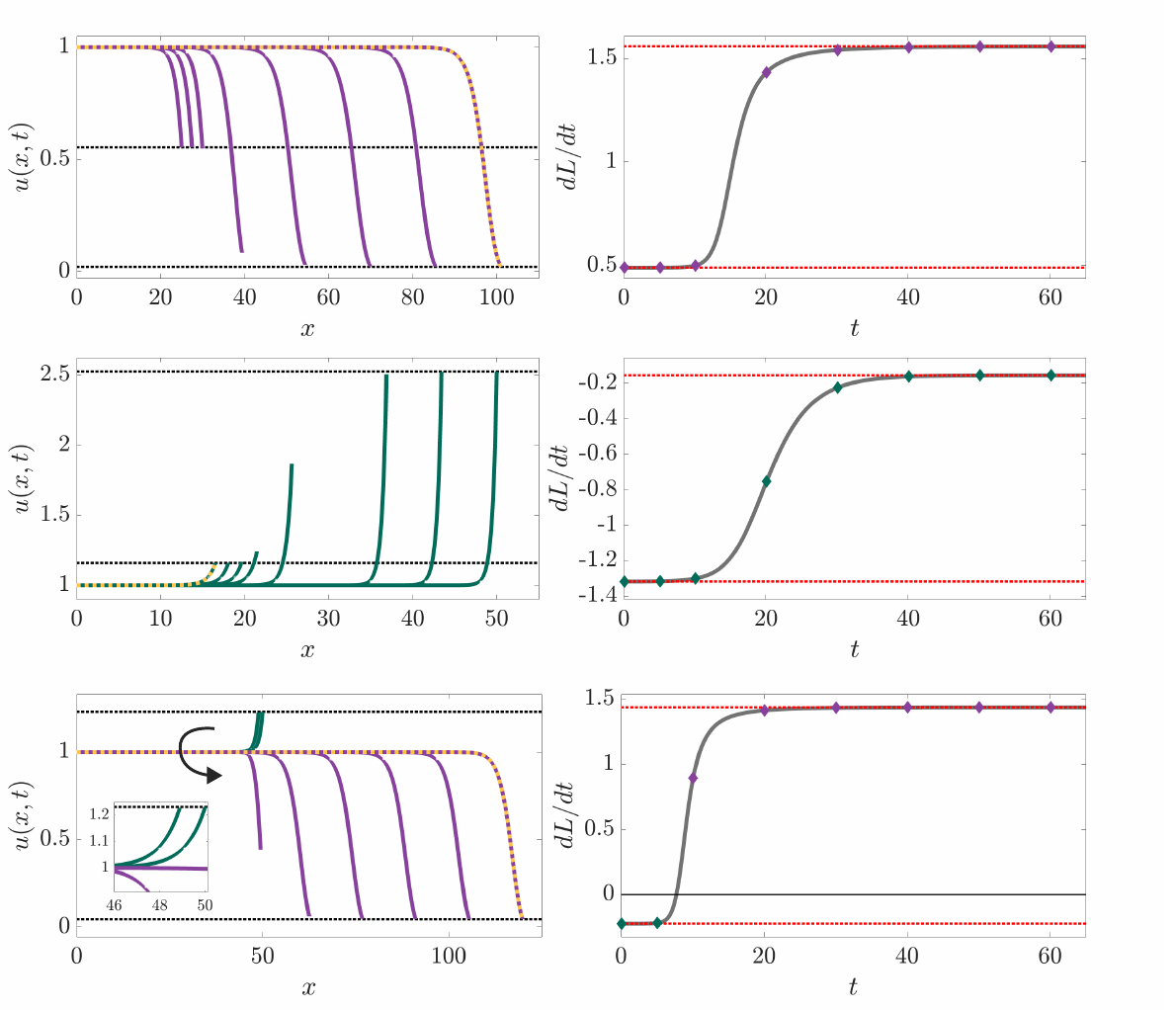}}%
    \put(0.01337381,0.82785254){\color[rgb]{0,0,0}\makebox(0,0)[t]{\lineheight{0.5}\smash{\begin{tabular}[t]{c}(a)\end{tabular}}}}%
    \put(0.01337381,0.56377078){\color[rgb]{0,0,0}\makebox(0,0)[t]{\lineheight{0.5}\smash{\begin{tabular}[t]{c}(c)\end{tabular}}}}%
    \put(0.48630267,0.82785254){\color[rgb]{0,0,0}\makebox(0,0)[t]{\lineheight{0.5}\smash{\begin{tabular}[t]{c}(b)\end{tabular}}}}%
    \put(0.48630267,0.56377078){\color[rgb]{0,0,0}\makebox(0,0)[t]{\lineheight{0.5}\smash{\begin{tabular}[t]{c}(d)\end{tabular}}}}%
    \put(0.01337381,0.28449988){\color[rgb]{0,0,0}\makebox(0,0)[t]{\lineheight{0.5}\smash{\begin{tabular}[t]{c}(e)\end{tabular}}}}%
    \put(0.48630267,0.28449988){\color[rgb]{0,0,0}\makebox(0,0)[t]{\lineheight{0.5}\smash{\begin{tabular}[t]{c}(f)\end{tabular}}}}%
    \put(0,0){\includegraphics[width=\unitlength,page=2]{XFigure_8.pdf}}%
  \end{picture}%
\endgroup%

    \caption{(a,c,e) Numerical solutions of the PDE system \eqref{gov PDE}–\eqref{BC L} for different values of $\alpha$ and $\beta$, shown at $t = [0, 5, 10, 20, 30, 40, 50, 60]$. Direction of travel is marked (arrow). All chosen $(\alpha, \beta)$ pairs have two possible travelling wave solutions (see Fig. \ref{fig: paramspace fig}). All initial conditions are the expected wave form of the solution with the smaller wave speed, obtained numerically from the phase plane. Superimposed are the second wave profiles from the numerical phase plane trajectory (dashed yellow line). The two possible front densities $u^*$ from the phase plane analysis are shown (dashed black lines). (b,d,f) Front speeds at $x = L(t)$, with time points from (a,c,e) indicated (coloured markers). The two possible wave speeds $c$ from the phase plane analysis are shown (dashed red lines). In (a,b), $\alpha = -2$, $\beta = 1.6$, with $u_{1}^* = 0.55$, $u_{2}^* = 0.02$  and $c_1 = 0.49$, $c_2 = 1.56$. In (c,d), $\alpha = -0.85$, $\beta = 0.83$, with $u_{1}^* = 2.53$, $u_{2}^* = 1.16$  and $c_1 = -1.32$, $c_2 = -0.16$. In (e,f), $\alpha = -1.4$, $\beta = 1.5$, with $u_{1}^* = 1.23$, $u_{2}^* = 0.04$  and $c_1 = -0.22$, $c_2 = 1.44$. Colours (purple for $c > 0$, green for $c < 0$) indicate the sign of $c$. Note that the simulations of (e,f) correspond to the parameters shown in Fig. \ref{fig: BC intersections}c.}
    \label{fig: PDE double TW}
\end{figure}

The only finite regions of Fig. \ref{fig: paramspace fig} are those with two possible travelling wave solutions, say $c = c_{1}$ and $c = c_{2}$, where without loss of generality $c_{1} < c_{2}$. We show numerical results in Fig. \ref{fig: PDE double TW} for three possible scenarios that arise when $\alpha < -1/\sqrt{2}$ and $\beta < 2$: two invading solutions with $0 < c_{1} < c_{2}$, two receding solutions with $c_{1} < c_{2} < 0$, or one invading solution and one receding solution with $c_{1} < 0 < c_{2}$.
To observe both solutions, we specify the initial condition as the wave form of the solution with the wave speed $c_{1}$ and front density $u_{1}^{*}$. While the solution stays on this travelling wave for some time, it eventually evolves to the stable travelling wave solution with wave speed $c_{2}$. The shape of the initial condition is determined by appropriately shifting the density profile $U(z)$ for the $c_{1}$ solution, obtained in the same manner as for Fig. \ref{fig: result fig 1}. 

In Fig. \ref{fig: PDE double TW}a,b, where $\alpha = 2$ and $\beta = 1.6$, the solution is initialised as the travelling wave solution with front density $u_{1}^* = 0.55$ and wave speed $c_{1} = 0.49$. Over time, the front density decreases and the wave speed increases before the solution stabilises to the second travelling wave solution with $u_{2}^* = 0.02$ and $c_{2} = 1.56$. Similar dynamics are observed in Fig. \ref{fig: PDE double TW}c,d where $\alpha = -0.85$ and $\beta = 0.83$ with two receding solutions. The solution is initialised as a receding travelling wave solution with front density $u_{1}^* = 2.53$ and wave speed $c_{1} = -1.32$. After some time, the front density decreases and the wave speed becomes less negative before the solution stabilises as the travelling wave solution with $u_{2}^* = 1.16$ and $c_{2} = -0.16$. In Fig. \ref{fig: PDE double TW}e,f for $\alpha = -1.4$ and $\beta = 1.5$, the initial condition of a receding solution with front density $u_{1}^* = 1.23$ and wave speed $c_{1} = -0.22$ evolves to an invading wave with $u_{2}^* = 0.04$ and $c_{2} = 1.44$. Across these simulations we have shown the $c_{1}$ solution to be unstable, and for some
$(\alpha,\beta)$ pairs, a consequence of this instability is a complete change in the direction of movement. 

\subsection{Non-travelling wave behaviours}
\begin{figure}[!h]
\def\svgwidth{\textwidth}
 %% Creator: Inkscape 1.3 (0e150ed, 2023-07-21), www.inkscape.org
%% Accompanies image file 'XFigure_9.pdf' (pdf, eps, ps)
\begingroup%
  \makeatletter%
  \providecommand\color[2][]{%
    \errmessage{(Inkscape) Color is used for the text in Inkscape, but the package 'color.sty' is not loaded}%
    \renewcommand\color[2][]{}%
  }%
  \providecommand\transparent[1]{%
    \errmessage{(Inkscape) Transparency is used (non-zero) for the text in Inkscape, but the package 'transparent.sty' is not loaded}%
    \renewcommand\transparent[1]{}%
  }%
  \providecommand\rotatebox[2]{#2}%
  \newcommand*\fsize{\dimexpr\f@size pt\relax}%
  \newcommand*\lineheight[1]{\fontsize{\fsize}{#1\fsize}\selectfont}%
  \ifx\svgwidth\undefined%
    \setlength{\unitlength}{562.64351432bp}%
    \ifx\svgscale\undefined%
      \relax%
    \else%
      \setlength{\unitlength}{\unitlength * \real{\svgscale}}%
    \fi%
  \else%
    \setlength{\unitlength}{\svgwidth}%
  \fi%
  \global\let\svgwidth\undefined%
  \global\let\svgscale\undefined%
  \makeatother%
  \begin{picture}(1,0.59689991)%
    \lineheight{1}%
    \setlength\tabcolsep{0pt}%
    \put(0,0){\includegraphics[width=\unitlength,page=1]{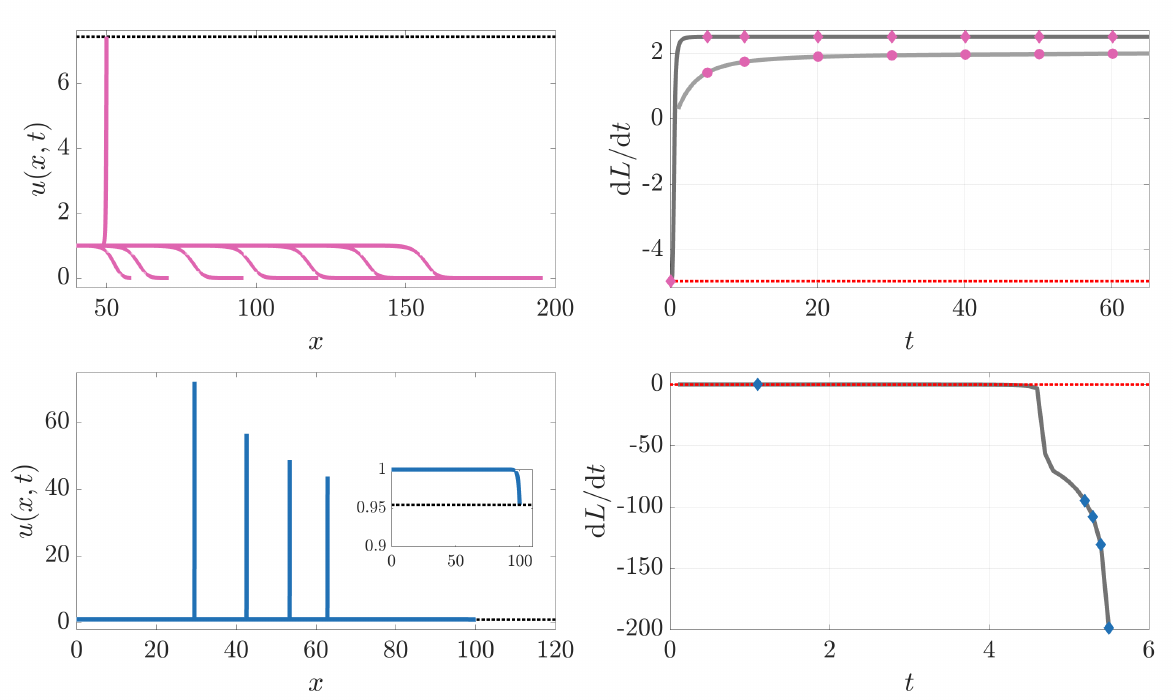}}%
    \put(0.01357074,0.56530207){\color[rgb]{0,0,0}\makebox(0,0)[t]{\lineheight{0.5}\smash{\begin{tabular}[t]{c}(a)\end{tabular}}}}%
    \put(0.01357074,0.29733158){\color[rgb]{0,0,0}\makebox(0,0)[t]{\lineheight{0.5}\smash{\begin{tabular}[t]{c}(c)\end{tabular}}}}%
    \put(0.52012349,0.56530207){\color[rgb]{0,0,0}\makebox(0,0)[t]{\lineheight{0.5}\smash{\begin{tabular}[t]{c}(b)\end{tabular}}}}%
    \put(0.52012349,0.29733158){\color[rgb]{0,0,0}\makebox(0,0)[t]{\lineheight{0.5}\smash{\begin{tabular}[t]{c}(d)\end{tabular}}}}%
    \put(0,0){\includegraphics[width=\unitlength,page=2]{XFigure_9.pdf}}%
  \end{picture}%
\endgroup%

    \caption{(a,c) Numerical solutions of the PDE system \eqref{gov PDE}–\eqref{BC L} for different values of $\alpha$ and $\beta$. Both of the chosen $(\alpha,\beta)$ pairs have a single travelling wave solution (see Fig. \ref{fig: paramspace fig}). All initial conditions are the expected wave profile, obtained numerically from the phase plane. Direction of travel is marked (arrow). The front density $u^*$ from phase plane analysis is shown (dashed black line). (b,d) Front speeds at $x = L(t)$, with time points from (a,c) indicated (diamond markers). In (b), we also show the speed (round markers) identified numerically by tracking the density $u = 0.5$ in (a). The wave speed $c$ from phase plane analysis is shown (dashed red line). In (a,b), $\alpha = -1$, $\beta = 2.5$, with $u^* = 7.45$ and $c = -4.95$. In (c,d), $\alpha = -2.1$, $\beta = 2.05$, with $u^* = 0.95$ and $c = 0.05$.}
    \label{fig: unstable PDE}
\end{figure}

We now show the non-travelling wave behaviours that emerge when there is no stable travelling wave solution. To do so, we consider the regions of Fig. \ref{fig: paramspace fig} where $\alpha < -1\sqrt{2}$ and $\beta > 2$ and single, receding or invading, travelling wave solutions exist. We specify the initial condition as the expected wave shape of the solution with wave speed $c$ and front density $u^*$. In Fig. \ref{fig: unstable PDE}a,b we set $\alpha = -1$ and $\beta = 2.5$, associated with a receding wave with speed $c = -4.95$ and front density $u^* = 7.45$. Rather than recede with front density $u^* \gg 1$, the front density rapidly decreases so that $u^* \to 0$. In the boundary velocity function \eqref{BC_fu} the influence of $\alpha$ becomes negligible, so that the boundary moves with speed $c = \beta$, where $\beta = 2.5$. As observed in Fig. \ref{fig: unstable PDE}a, the bulk body of the solution moves as a Fisher-KPP type travelling wave. Tracking the density $u = 0.5$, this bulk body invades with speed approaching $c_{0.5} = 2$. As $c > c_{0.5}$ (see Fig. \ref{fig: unstable PDE}b), the population front at $L(t)$ invades at a faster rate than the bulk body of the wave, increasing the distance between the two over time. 

In Fig. \ref{fig: unstable PDE}c,d we choose $\alpha = -2.1$ and $\beta = 2.05$, associated with an unstable, invading travelling wave with $c = 0.05$ and $u^* = 0.95$. Instead of establishing an invading wave, the front density $u^*$ rapidly increases above the carrying capacity so that $u^* \gg 1$. This causes $\alpha$ to become the dominant influence of the boundary velocity. Given $\alpha < 0$, the solution recedes, with the mass-conserving boundary condition at $x = L(t)$ increasing the front density $u^*$. This further amplifies the boundary velocity, establishing a positive feedback loop whereby the solution rapidly recedes towards the boundary at $x = 0$. 

\section{Stability of solutions along $\alpha = -\beta$}\label{sec: stability analysis}

\begin{figure}[!h]
\def\svgwidth{\textwidth}
 \input{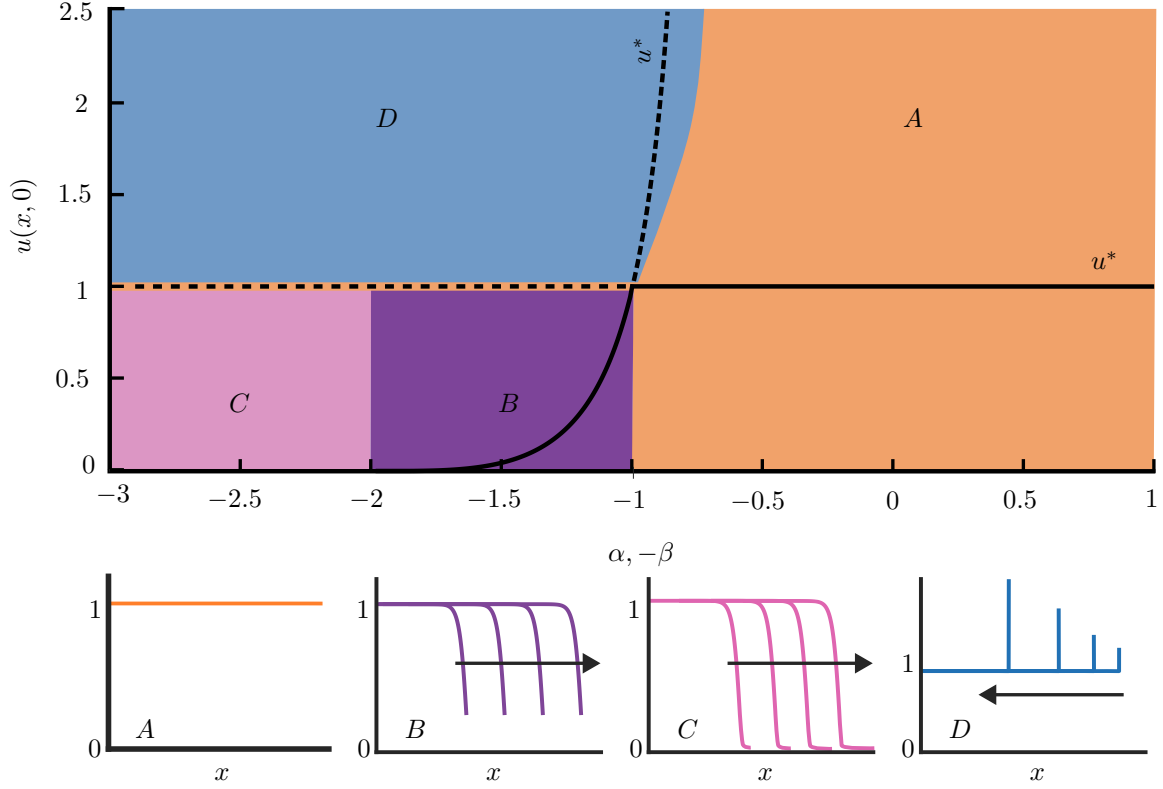}
    \caption{Numerical analysis of steady state stability along $\alpha = -\beta$ for various uniform initial conditions. On the horizontal axis we vary $\alpha$ (and $-\beta$), and on the vertical axis we vary the uniform initial density $u(x,0) = u_{0}$. The uniform initial conditions that we consider evolve to four different behaviours, shown in the corresponding schematics and described as follows. Behaviour A is the steady state solution and it is stable for $\alpha > -1/\sqrt{2}$. Behaviour B is the invading travelling wave solution, possible from $-2 < \alpha < -1$ (recall Fig. \ref{fig: paramspace fig}). In Behaviour C, where $\beta > 2$ ($\alpha < -2$), the invasive front extends with speed $c = \beta$ while the bulk wave speed approaches $c = 2$ (observed in Fig. \ref{fig: unstable PDE}a,b). Finally, in Behaviour D the front density increases and the front speed becomes more negative in a positive feedback loop (observed in Fig. \ref{fig: unstable PDE}c,d). We also superimpose the steady state solution $u^* = 1$, and invading $(u^* < 1)$ and receding $(u^* > 1)$ travelling wave solutions (black lines, solid when stable and dashed when unstable). We observe a transcritical bifurcation, where the stability of the steady state solution exchanges stability with the invading travelling wave solution at $\alpha = -1$. For $-1 < \alpha < -1/\sqrt{2}$, we observe that the behaviour of the system is driven by the value of the uniform initial density $u_{0}$ relative to the basin boundary (undetermined) separating the steady state and travelling wave solutions. We note that since the initial condition under consideration is uniform, this boundary is not expected to exactly coincide with the value of $u^*$ of the non-uniform initial condition. As $\alpha$ approaches $-2$, the basin of attraction of the steady state solution decreases so that behaviour $D$, associated with the unstable travelling wave solution, is increasingly favoured.}
    \label{fig: steady state stability}
\end{figure}

Observations of numerical simulations of the PDE system \eqref{gov PDE}–\eqref{BC L} described in the previous section have shown the existence of unstable travelling wave solutions. Unstable solutions first emerged when we considered $\alpha < -1/\sqrt{2}$ and $\beta < 2$ in Fig. \ref{fig: PDE double TW} to observe the coexistence of stable and unstable travelling wave solutions. With the additional constraint that $\beta > 2$ in Fig. \ref{fig: unstable PDE}, we could no longer identify stable travelling wave solutions at all. A complete analysis of the stability of invariant solutions in $(\alpha, \beta)$ space is beyond the scope of this work. However, the stability of the system along $\alpha = -\beta$ can be studied more readily, since this relationship is satisfied by the steady state solution $u(x,t) = 1$. We focus on this closed form solution as determining the linear stability of other travelling waves would require numerical methods. 

We first observe the stability of this uniform state by setting $\alpha = -\beta$ in the boundary velocity function and varying the initial condition, $u(x,0) = u_{0}$, where $u_{0}$ is a constant. The numerical behaviours are summarised in Fig. \ref{fig: steady state stability}. All initial conditions that we consider tend to the stable steady state for $\alpha > -1\sqrt{2}$. Initial conditions sufficiently close to $u(x,t) = 1$ remain stable for $-1 < \alpha < -1/\sqrt{2}$, and develop into the receding runaway solution that was observed in Fig. \ref{fig: unstable PDE}c,d for large $u_{0}$. This is consistent with the presence of an unstable travelling wave that separates attraction to the steady state from the runaway behaviour. In this regime, the speed of the unstable receding travelling wave increases in magnitude as $\alpha \to -1/\sqrt{2}$. As $\alpha \to -1$, this receding travelling wave slows and approaches zero speed, consistent with a collision with the steady state in a transcritical bifurcation. For $\alpha < -1$, the stability of the steady state is lost to any perturbations of the initial condition $u(x,0) = 1$.  For $u_{0} < 1$ we observe an invading travelling wave solution, which exists between $-2 <\alpha < -1$ (see Fig. \ref{fig: paramspace fig}), and for $u_{0} > 1$ we observe a receding runaway solution. Finally, for $\alpha < -2$, which corresponds to $\beta > 2$, there is no stable invading solution, so that for $u_{0} < 1$ the system takes on the behaviour shown in Fig. \ref{fig: unstable PDE}a,b, where the population front approaches zero and invades at a higher speed than the bulk body of the population. 

Our numerical analysis of the stability for $\alpha = -\beta$ is supported analytically by considering the steady state solution with a small perturbation 
\begin{align}
    u(x,t) = 1 + \tilde{u}(x,t).\label{pert steady state}
\end{align}
The leading order front $L(t)$ is set to be zero without loss of generality. From substituting (\ref{pert steady state}) into the system \eqref{gov PDE}–\eqref{BC L} and taking linear terms only, the problem for $\tilde u$ is
\begin{align}
    \tilde{u}_{t} = \tilde{u}_{xx} - \tilde{u}, \quad -\infty < x < 0,\label{SS gov equ}
\end{align}
with $\tilde{u} \to 0$ as $x \to -\infty$ and at $x = 0$, 
\begin{align}
    \tilde{u}_{x}(0,t) + \alpha\tilde{u}(0,t) = 0. \label{SS BC}
\end{align}
This is a linear reaction-diffusion problem (\ref{SS gov equ}) on a semi-infinite domain with a Robin boundary condition (\ref{SS BC}).  When $\alpha > 0$, the solution may be found using a Fourier transform method \cite{pinsky2011partial}. We define a generalised Fourier transform which takes into account the boundary condition (\ref{SS BC}):
\begin{align}
    \hat{u}(k,t) &= \mathcal{F}_{\alpha}\{{\tilde{u}(x,t)}\}= \int_{-\infty}^{0} {\tilde{u}(x,t)}(k\textrm{cos}(kx) - \alpha\textrm{sin}(kx))\textrm{d}x.
\end{align}
Applying this transform to the governing equation (\ref{SS gov equ}) results in a linear ODE for the transformed variable $\hat u$:
\begin{align*}
    \hat{u}_{t} = -(k^2 + 1)\hat{u},
\end{align*}
which has the solution 
\begin{align*}
\hat{u}(k,t) = \hat{u}(k,0)e^{-(k^2 + 1)t}.
\end{align*}
Here $\hat u(k,0) = \mathcal F_\alpha\{u(x,0)\}$ is the transform of the initial condition. Inverting the transform, we thus obtain the solution 
\begin{align*}
    \tilde{u}(x,t) = \mathcal F_\alpha^{-1}\{\hat u(k,t)\} = \int_{0}^{\infty}\frac{2}{\pi(\alpha^2 + k^2)}\hat{u}(k,0)e^{-(k^2+1)t}(k\cos{(kx)} - \alpha \sin{(kx)})\mathrm{d}k.
\end{align*}
which will decay as it is bounded by $e^{-t}\tilde{u}(x,0)$.

While this solution is valid for $\alpha > 0$, we are interested in the case where $\alpha < 0$, for which the PDE has a particular solution 
\begin{align}
    \tilde{u}_{1}(x,t) = Ce^{-\alpha x + (\alpha^2 - 1)t}. \label{SS particular}
\end{align}
This solution is orthogonal to the Fourier integral solution, that is, $\mathcal F_a(\tilde u_1(x,t)) = 0$.  This extra term must therefore be included in the general solution when $\alpha < 0$. The general solution for $\alpha < 0$ is then 
\begin{align}
    \tilde{u}(x,t) = \int_{0}^{\infty}\frac{2}{\pi(\alpha^2 + k^2)}\hat{u}(k,0)e^{-(k^2+1)t}(k\cos{kx} - \alpha \sin{kx{)\mathrm{d}k}} + Ce^{-\alpha x + (\alpha^2 - 1)t}. \label{SS final eq}
\end{align}
Here the constant $C$ is related to the initial condition by
\begin{align}
C = 2\alpha \int_{-\infty}^0 e^{-\alpha x} \tilde u(x,0) \, \mathrm dx.
\end{align}
The contribution of the final term of (\ref{SS final eq}), arising from the particular solution (\ref{SS particular}), decays for $-1 < \alpha < 0$ and grows for $\alpha <-1$. This confirms our numerical finding that when $\alpha < -1$, the steady state is not stable for any small perturbation to the exact solution. While this analysis is focused on the relationship $\alpha = -\beta$, we find that it is representative of the stability properties that we observed more generally from numerical simulations presented in Section \ref{sec: PDE sims}.

\section{Discussion}

In this article, we introduce a mass-conserving free boundary condition to the Fisher–KPP equation. In contrast to the Fisher–Stefan model, in which boundary motion is governed by a flux condition involving a leakage coefficient $\kappa$, our formulation enforces conservation of mass at the moving boundary and prescribes the boundary velocity as a function of the local cell density. This mass-conserving condition is motivated by settings in which cells remodel their environment without being consumed, and arises naturally in existing models of biological phenomena \cite{gaffney1999modelling,chen2000free,baker2019free,murphy2021travelling}. By exploiting the well-established travelling wave structure of the Fisher-KPP equation, we readily examined the influence of the mass-conserving boundary condition on the resulting population dynamics to identify the existence of single, and multiple travelling wave solutions.

A key distinction of our model is that satisfying the mass-conserving boundary condition in the phase plane requires trajectory truncations $(u^*,-cu^*)$ that are both wave speed and front density dependent. We observed travelling wave solutions characterised by both $-\infty < c < 2$ and $0 < u^* < \infty$, where there is a one-to-one relationship between the two (see wave speed curve in Fig. \ref{fig: result fig 2}). Invading solutions have a front density below carrying capacity $u^* \in (0,1)$ while receding solution have a front density above carrying capacity $u^* \in (1,\infty).$ This contrasts with the Fisher-Stefan model, where the boundary condition is prescribed solely by the wave speed $c$ and the front density is always $u^* = 0$. Our formulation may be advantageous in settings where the boundary dynamics cannot be observed directly. As discussed in  El-Hachem et al. \cite{El-Hachem_McCue_Simpson_2021_I_R}, estimates of a population's wave speed $c$ are readily available, given they can be measured experimentally \cite{Maini_McElwain_Leavesley_2004,Simpson_Zhang_Mariani_Landman_Newgreen_2007}. We could use experimental observations of either the wave speed or front density to determine a unique position $(u^*,c)$ along the wave speed curve, and the parameterisation of the boundary velocity function is constrained only by the requirement of passing through $(u^*,c)$. 

Utilising the wave speed curve we developed, we observed complex dynamics arising from a simple linear boundary velocity function. These include single invading and receding solutions, as well as multiple travelling wave solutions, of which only the solution with the most positive wave speed is stable. In the parameter regime with multiple solutions, it is inherently easier to study stable invasive behaviours. Despite observing receding behaviour in early time, the direction of movement can reverse to converge to a stable invading solution. Interestingly, in the absence of a stable travelling wave solution, the system will either runaway towards the fixed boundary, or the population front will travel with speed $\beta$ to extend beyond the bulk body whose speed approaches the minimum nondimensional Fisher-KPP speed of $c = 2$. 

Through a combination of numerical and analytic techniques, we have developed a preliminary understanding of the stability properties of travelling wave solutions. A full stability analysis of the system is beyond the scope of this work. Instead, we focused on the steady state solution, $c = 0$ and $u = 1$ for $\alpha = -\beta$. The speed of an unstable receding solution approaches zero as $\alpha \to -1$ from above, and becomes a stable travelling wave for $\alpha < -1$. At $\alpha = -1$ this coincides with the steady state solution, resulting in a transcritical bifurcation. More generally across the $\alpha$, $\beta$ space, unstable travelling wave solutions arise if $\alpha < -1/\sqrt{2}$ or $\beta>2$. However, this does not impact our ability to simulate a stable travelling wave with arbitrary speed $-\infty < c < 2$ and front density $0 < u^* < \infty$. Imposing a stability condition $\beta < 2$, the boundary velocity function $c = \alpha u^* + \beta$ implies that $\alpha u^* > c- 2$. Since $c < 2$, and $u^* > 0$, this is always satisfied by $\alpha > 0$, ensuring a parameterisation will exist within a stable regime. 

A key advantage of this work is that the phase plane analysis we have undertaken is independent of the choice of boundary velocity function. Intuition surrounding travelling wave dynamics of the system under new boundary velocity functions can be developed by identifying intersections between the boundary function and the wave speed curve we have identified, as we demonstrated in Fig. \ref{fig: BC intersections}. A natural extension is to consider the boundary velocity to be a nonlinear function of cell density, using higher powers of $u$ or Hill functions. The framework we have developed is readily extendable to these boundary functions, offering a means of studying increasingly complex and biologically relevant boundary dynamics. Furthermore, we could extend the model to higher dimensions to consider the contribution of geometric effects \cite{tam2022effect,tam2023pattern}. 

\section*{Acknowledgements}

ALJ and GRW acknowledge funding from the Australian Research Council (ARC)
Discovery Early Career Researcher Award (DECRA) DE240100650 and the QUT
Faculty of Science Computational Bioimaging Group (CBG). ALJ also acknowledges the ARC Discovery Project (DP) DP230100025. MCD acknowledges support from the ARC Discovery Project DP250101095.

\bibliographystyle{plainnat}
\bibliography{references}

\end{document}